\documentclass[11pt]{article}

\usepackage[utf8]{inputenc}
\usepackage{amssymb,epsfig}
\usepackage{amsmath}
\usepackage{comment}
\usepackage[dvips]{color}
\usepackage[T1]{fontenc}

\newtheorem{defi}{Definition}
\newtheorem{prop}[defi]{Proposition}
\newtheorem{theo}[defi]{Theorem}
\newtheorem{conj}[defi]{Conjecture}
\newtheorem{lemm}[defi]{Lemma}
\newtheorem{coro}[defi]{Corollary}
\newtheorem{rema}[defi]{Remark}
\newtheorem{exem}[defi]{Example}
\newtheorem{exems}[defi]{Examples}

\newcommand{\bdefi}{\begin{defi}}
\newcommand{\edefi}{\end{defi}}
\newcommand{\bprop}{\begin{prop}}
\newcommand{\eprop}{\end{prop}}
\newcommand{\btheo}{\begin{theo}}
\newcommand{\etheo}{\end{theo}}
\newcommand{\blemm}{\begin{lemm}}
\newcommand{\brema}{\begin{rema}}
\newcommand{\erema}{\end{rema}}
\newcommand{\bexer}{\begin{exem}}
\newcommand{\eexer}{\end{exem}}
\newcommand{\bexems}{\begin{exems}}
\newcommand{\eexems}{\end{exems}}
\newcommand{\bconj}{\begin{conj}}
\newcommand{\econj}{\end{conj}}
\newcommand{\elemm}{\end{lemm}}
\newcommand{\bcoro}{\begin{coro}}
\newcommand{\ecoro}{\end{coro}}
\newcommand{\dem}{\noindent{\bf Proof. }}
\newcommand{\rem}{\noindent{\bf Remark. }}

\usepackage{mathrsfs}
\renewcommand\mathcal{\mathscr}

\newcommand{\B}{{\cal B}}
\renewcommand{\L}{{\cal L}}
\newcommand{\M}{{\cal M}}

\newcommand{\OOO}{{\cal O}}
\newcommand{\C}{{\cal C}}

\newcommand{\maths}[1]{{\mathbb #1}}  

\newcommand{\RR}{\maths{R}}
\newcommand{\NN}{\maths{N}}
\newcommand{\CC}{\maths{C}}
\newcommand{\QQ}{\maths{Q}}

\newcommand{\KK}{\maths{K}}
\newcommand{\ZZ}{\maths{Z}}

\newcommand{\PP}{\maths{P}}

\newcommand{\LL}{\maths{L}}

\newcommand{\aaa}{{\mathfrak a}}
\renewcommand{\ggg}{{\mathfrak g}}

\newcommand{\qqq}{{\mathfrak q}}
\renewcommand{\lll}{{\mathfrak l}}

\newcommand{\uuu}{{\mathfrak u}}

\newcommand{\ppp}{{\mathfrak p}}
\newcommand{\hhh}{{\mathfrak h}}

\newcommand{\mmm}{{\mathfrak m}}
\renewcommand{\uuu}{{\mathfrak u}}

\newcommand{\UUU}{{\mathfrak U}}

\newcommand{\ra}{\rightarrow}
\newcommand{\bs}{\backslash}

\newcommand{\wt}[1]{{\widetilde{#1}}}

\newcommand{\ga}{\gamma}
\newcommand{\Ga}{\Gamma}

\newcommand\vre{\varepsilon}

\newcommand{\cqfd}{\hfill$\Box$}

\newcommand{\card}{{\operatorname{Card}}}

\newcommand{\Vol}{\operatorname{Vol}}

\newcommand{\SL}{\operatorname{SL}}
\newcommand{\GL}{\operatorname{GL}}

\newcommand\Ad{\operatorname{Ad}}

\newcommand\vol{\operatorname{vol}}

\newcommand\GalQ{\operatorname{Gal}(\overline{\QQ}/\QQ)}

\newcommand{\ellk}{{q}}

\newcounter{fig}

\title{Counting orbits of integral points \\
in families of affine homogeneous varieties\\ and diagonal flows}
\author{Alexander Gorodnik \and Fr\'ed\'eric Paulin} 
\date{\today}

\begin{document}
\bibliographystyle{../alphanum}

\maketitle

\begin{abstract}
  In this paper, we study the distribution of integral points on
  parametric families of affine homogeneous varieties. By the work of
  Borel and Harish-Chandra, the set of integral points on each such
  variety consists of finitely many orbits of arithmetic groups, and
  we establish an asymptotic formula (on average) for the number of
  the orbits indexed by their Siegel weights. Our arguments use the
  exponential mixing property of diagonal flows on homogeneous spaces.
\end{abstract}

\section{Introduction}
\label{sect:intro}

In this paper, we investigate the asymptotic distribution of integral
points on families of homogeneous algebraic varieties using dynamical
systems techniques.  Let $\bf L$ be a reductive algebraic group in
$\GL_n(\CC)$ defined over $\QQ$, and let ${\bf X}_v={\bf L}v$ with
$v\in \mathbb{Q}^n$ be a Zariski-closed set. Borel-Harish-Chandra's
finiteness theorem \cite[Theo.~6.9]{BorelHarishChandra62} says that
the number of orbits of ${\bf L}(\ZZ)$ in ${\bf X}_v(\ZZ)$ is finite.
We investigate the asymptotic behaviour of this number where the
orbits are counted with suitable weights. Namely, to each orbit
$[u]\in {\bf L}(\ZZ)\bs{\bf X}_v(\ZZ)$, we associate the Siegel weight
$$
w([u])=\frac{\vol\big( \hbox{Stab}_{\bf L}(u)(\ZZ)\bs
\hbox{Stab}_{\bf L}(u)(\RR)\big)}
{\vol\big( {\bf L}(\ZZ)\bs
{\bf L}(\RR)\big)}
$$
and consider
$$
\beta(v)=\sum_{[u]\in 
{\bf L}(\ZZ)\bs{\bf X}_v(\ZZ)} w([u]).
$$
This quantity appears naturally in the study of spatial distribution
of integral points.  Given a cone $\Omega\subset \RR^n$ and a sequence
of colinear vectors $v_n \in \QQ^n$ such that ${\bf X}_{v_n} (\RR)
\cap \Omega$ is bounded with boundary of measure zero and
$\hbox{Stab}_{\bf L}(v_n)$ is a maximal subgroup of $\bf L$, it was
shown in \cite{EskOh06}, when $\bf L$ is semi-simple, that
\begin{equation}\label{eq:aaaaa}
|{\bf X}_{v_n}(\mathbb{Z})\cap \Omega|\sim \vol(\Omega)\cdot \beta(v_n)
\quad\hbox{as $\beta(v_n)\to \infty$.}
\end{equation}
Under some additional conditions on the varieties ${\bf X}_{v_n}$, the
asymptotics of $|{\bf X}_{v_n}(\mathbb{Z})\cap \Omega|$ can be also
computed using the Hardy--Littlewood circle method, so that the
quantities $\beta(v_n)$ can be expressed in terms of local densities
(see \cite{BorRud95, Oh04}).  However, we emphasise that the
asymptotic formula (\ref{eq:aaaaa}) and the Hardy--Littlewood method
do not apply in general. In particular, it might happen that
$\beta(v_n)$ does not converge to $\infty$ even though $\beta(v_n)$
does on average. Contrarily to the last two references, we also
consider non semi-simple cases, which are useful for some examples.

The aim of this paper is to develop a direct argument that establishes
an asymptotic formula for the quantities $\beta(v)$ on average.  Our
methods use flows on homogeneous spaces and dynamical systems
techniques.  We show that sums over $\beta(v)$ can be interpreted as
volumes of intersections of two transversal submanifolds in a suitable
homogeneous space. Then we prove that one of these submanifolds
becomes equidistributed and deduce an asymptotic formula for the sums
of $\beta(v)$'s.

\medskip Now we proceed to state our result precisely.  Let $\bf L$ be
a reductive linear algebraic group defined and anisotropic over $\QQ$,
let $\pi:{\bf L}\ra \GL({\bf V})$ be a rational linear representation
of $\bf L$ defined over $\QQ$ and let $\Lambda$ be a $\ZZ$-lattice in
${\bf V}(\QQ)$ invariant under ${\bf L}(\ZZ)$. For every $v\in {\bf
  V}(\QQ)$ whose orbit ${\bf X}_v$ under $\bf L$ is Zariski-closed in
$\bf V$, the number of orbits of ${\bf L}(\ZZ)$ in ${\bf X}_v\cap
\Lambda$ is finite. We will count these orbits using appropriate
weights. For every $u\in {\bf X}_v(\QQ)$ with stabiliser ${\bf L}_u$
in $\bf L$, define the {\it Siegel weight} of $u$ as
$$
w_{{\bf L},\pi}(u)=
\frac{\vol\big( {\bf L}_u(\ZZ)\bs
{\bf L}_u(\RR)\big)}
{\vol\big( {\bf L}(\ZZ)\bs
{\bf L}(\RR)\big)}\;,
$$
using Weil's convention for the normalisation of the measures on ${\bf
  L}_u(\RR)$ (depending on the choice of a left Haar measure on ${\bf
  L}(\RR)$ and of a ${\bf L}(\RR)$-invariant measure on ${\bf
  X}_v(\RR)$, see Section \ref{sect:mainproof}). These weights
generalise the ones occuring in Siegel's weight formula when $\bf L$
is an orthogonal group (see for instance \cite{Siegel44,EskRudSar91},
and \cite[Chap.~5]{Voskresenskii98} for general $\bf L$).

A particular case of the main results of this paper is the following
one.
\btheo \label{theo:intromaincourtsc} Let $\bf G$ be a simply connected
reductive linear algebraic group defined over $\QQ$, without
nontrivial $\QQ$-characters.  Let $\bf P$ be a maximal parabolic
subgroup of $\bf G$ defined over $\QQ$, and let ${\bf P} ={\bf AMU}$
be a relative Langlands decomposition of $\bf P$, such that ${\bf A}
(\RR)_0$ is a one-parameter subgroup $(a_s)_{s\in\RR}$, with
$\lambda=\log \det\, (\Ad a_1)_{\mid \UUU}>0$, where $\UUU$ is the Lie
algebra of ${\bf U} (\RR)$.  Let $\rho:{\bf G} \ra \GL({\bf V})$ be a
rational representation of $\bf G$ defined over $\QQ$ such that there
exists $v_0\in{\bf V}(\QQ)$ whose stabiliser in $\bf G$ is $\bf
MU$. Let $\bf L$ be a reductive algebraic subgroup of $\bf G$ defined
and anisotropic over $\QQ$. Assume that $\bf LP$ is Zariski-open in
$\bf G$ and that for every $s\in\RR$, the orbit ${\bf X}_s=\rho({\bf
  L} a_s)v_0$ is Zariski-closed in $\bf V$.

Let $\Lambda$ be a $\ZZ$-lattice in ${\bf V}(\QQ)$ invariant under
${\bf G}(\ZZ)$, and let $\Lambda^{\rm prim}$ be the subset of
indivisible elements of $\Lambda$. Assume $\rho$ to be irreducible
over $\CC$. Then there exist $c,\delta>0$ such that, as $t$ tends to
$+\infty$,
$$
\sum_{0\leq s\leq t}\;\;\;\sum_{[x]\in 
{\bf L}(\ZZ)\bs({\bf X}_s\cap \Lambda^{\rm prim})} \;
w_{{\bf L},\rho_{\mid {\bf L}}}(x)= c\;e^{\lambda t}+
\operatorname{O}(e^{(\lambda-\delta)t})\;.
$$
\etheo

More precisely, let ${\bf G,P,A,M,U,V,L}, \rho, v_0, (a_s)_{s\in\RR}$
be as above ($\rho$ not necessarily irreducible). Endow ${\bf G}
(\RR)$ with a left-invariant Riemannian metric, for which the Lie
algebras of ${\bf MU}(\RR)$ and ${\bf A}(\RR)$ are orthogonal, and the
orthogonal of the Lie algebra of ${\bf P}(\RR)$ is contained in the
Lie algebra of ${\bf L}(\RR)$.

\btheo \label{theo:intromainsc} There exists $\delta>0$ such that,
as $t$ tends to $+\infty$,
\begin{align*}
\sum_{0\leq s\leq t}\;\;\;&\sum_{[x]\in 
{\bf L}(\ZZ)\bs(\rho({\bf L}(\RR) a_s)v_0\cap \rho({\bf G}(\ZZ))v_0)} \;
w_{{\bf L},\rho_{\mid {\bf L}}}(x)\\ &\;\;\;\;\;\;\;\;=\;\; 
\frac{\vol\big({\bf MU}(\ZZ)\bs {\bf MU}(\RR)\big)
\vol(a_1^\ZZ\bs {\bf A}(\RR)_0)}
{\lambda\vol({\bf G}(\ZZ)\bs {\bf G}(\RR)}\;e^{\lambda t}+
\operatorname{O}(e^{(\lambda-\delta)t})\;.
\end{align*}
\etheo

We will prove a more general version of this result in Section
\ref{sect:mainproof} without the maximality condition on $\bf P$,
involving the more elaborate root data of $\bf P$, and without the
simple connectedness assumption on $\bf G$ (up to a slight
modification of the Siegel weights), see Theorem \ref{theo:main} and
Theorem \ref{theo:mainmaximal}.

\medskip We illustrate Theorem \ref{theo:intromaincourtsc} by an
example, which is new: we give an asymptotic estimate on the
(weighted) number of inequivalent integral points on hyperplane
sections of affine quadratic surfaces.  More examples are given in
Section \ref{sect:applications}.

\bcoro \label{theo:introquadraticsurface} Let $n\geq 3$, let
$q:\CC^n\ra \CC$ be a nondegenerate rational quadratic form, which is
isotropic over $\QQ$, let $\ell:\CC^n\ra \CC$ be a nonzero rational
linear form, and let ${\bf L}= \{g\in\SL_n(\CC) \;:\;q\circ
g=q,\;\ell\circ g=\ell \}$.  For every $k\in\QQ$, let $\Sigma_k$ be the
set of primitive $x\in\ZZ^n$ such that $q(x)=0$ and
$\ell(x)=k$. Assume that the restriction of $q$ to the kernel of
$\ell$ is nondegenerate and anisotropic over $\QQ$. Then there exist
$c=c(q,\ell)>0$ and $\delta=\delta(q)>0$ such that, as $r\ra+\infty$,
$$
\sum_{k\in[1,r]} \;\;\;\sum_{[u]\in 
({\bf L}(\ZZ)\cap {\bf L}(\RR)_0)\bs\Sigma_k}\;
\vol\big( ({\bf L}_u(\ZZ)\cap{\bf L}(\RR)_0)\bs
({\bf L}_u\cap{\bf L}(\RR)_0)\big) = c\;r^{n-2}+
\operatorname{O}\big(r^{n-2-\delta}\big)\;.
$$
\ecoro

This result fits into the above program (up to a slight
modification of the Siegel weights, see Section \ref{sect:mainproof}),
by taking ${\bf V} =\CC^n$, $\Lambda=\ZZ^n$, and $\pi:{\bf L}\ra
\GL({\bf V})$ the inclusion map, noting that ${\bf L}$ is semisimple,
and defined and anisotropic over $\QQ$ as a consequence of the
assumptions (see Section \ref{sect:quadratic} for details, where we
also explicit $c$).

\medskip Asymptotics of the number of integral points in affine
homogeneous varieties has been extensively studied over the last
decades using harmonic analysis and dynamical systems techniques.  See
for instance \cite{EskRudSar91, DukRudSar93, EskMcMul93, EskMozSha96,
  GanOh03, Oh04, EskOh06}, as well as the surveys
\cite{Babillot02a,Oh10}.  Our results are quite different, since we
are counting whole orbits, weighted by the Siegel weights, of integral
points.  The asymptotic formula that we obtain here is similar to the
classical asymptotic formula for the number of integral quadratic
forms averaged over discriminant, proved by Siegel in
\cite{Siegel44b}.  More generally, if ${\bf V}$ is a prehomogeneous
vector space, analogous asymptotic formulas can be deduced from the
analytic properties of the corresponding zeta functions that were
studied by Sato and Shintani in \cite{SatShi74} (see also the
monograph \cite{Kimura03}).  However, the dimension of ${\bf G}v_0$ in
Theorem \ref{theo:intromaincourtsc} is typically much smaller than the
dimension of ${\bf V}$, and these methods don't apply.
We also
note that we do not assume ${\bf X}_v(\RR)$ to be an affine symmetric
space or that the stabiliser is a maximal subgroup, 
contrarily to \cite{DukRudSar93} and many other references.
Another difference with the counting results of
\cite{EskMozSha96,Oh04,EskOh06} is that these papers are using the
dynamics of unipotent flows, as instead we are using here the mixing
property with exponential decay of correlations of diagonalisable
flows, in the spirit of \cite{KleMar96} (see also
\cite{EskMcMul93,BenOh12}).  We are using the proof of the main result
of \cite{ParPau12JMD} as a guideline.

\medskip {\small {\it Acknowledgment. } We thank
  Jean-Louis Colliot-Thélène for discussions of our decomposable form
  result. The second author thanks the University of Bristol for a
  very fruitful short stay where this paper was conceived. The first
  author thanks the University Paris-Sud 11 for a month of invited
  professor where this paper was completed.}

\section{Counting Siegel weights}
\label{sect:mainproof}

Here are a few notational conventions. By linear algebraic group ${\bf
  G}'$ defined over a subfield $k$ of $\CC$, we mean a subgroup of
$\operatorname{GL}_N(\CC)$ for some $N\in\NN$ which is a closed
algebraic subset of $\M_N(\CC)$ defined over $k$, and we
define ${\bf G}' (\ZZ) ={\bf G}'\cap\operatorname{GL}_N(\ZZ)$. For
every linear algebraic group ${\bf G}'$ defined over $\RR$, we denote
by ${\bf G}' (\RR)_0$ the identity component of the Lie group of real
points of ${\bf G}'$. We denote by $\log$ the natural logarithm.

\medskip Let us first recall {\it Weil's normalisation} of measures on
homogeneous spaces. Let $G'$ be a unimodular real Lie group, endowed
with a transitive smooth left action of $G'$ on a smooth manifold
$X'$, with unimodular stabilisers. A triple $(\nu_{G'},\nu_{X'},
(\nu_{G'_x})_{x\in X'})$ of a left Haar measure $\nu_{G'}$ on $G'$, a
left-invariant (Borel, positive, regular) measure $\nu_{X'}$ on $X'$
and of a left Haar measure $\nu_{G'_x}$ on the stabiliser $G'_x$ of
every $x\in X'$, is {\it compatible} if, for every $x\in X'$, for
every $f:G'\ra\RR$ continuous with compact support, with
$f_x:X'\ra\RR$ the map (well) defined by $gx\mapsto \int_{h\in
  G'_x}f(gh)\;d\nu_{G'_x}(h)$ for every $g\in G$, we have
$$
\int_{G'} f\;d\nu_{G'}=\int_{X'}f_x\;d\nu_{X'}\;.
$$
Weil proved (see for instance \cite[\S 9]{Weil65}) that, for every
left-invariant measure $\nu_{X'}$ on $X'$, then

$\bullet$ for every left Haar measure $\nu_{G'}$ on $G'$, there exists
a unique compatible triple $(\nu_{G'},\nu_{X'},(\nu_{G'_x})_{x\in X'})$.

$\bullet$ for every $x_0\in X'$, for every left Haar measure $\nu_0$ on
$G'_{x_0}$, there exists a unique compatible triple
$(\nu_{G'},\nu_{X'},(\nu_{G'_x})_{x\in X'})$ with
$\nu_{G'_{x_0}}=\nu_0$.

The following remark should be well-known, though we did not found a
precise reference.

\blemm\label{lem:weilnormal} If $(\nu_{G'},\nu_{X'},(\nu_{G'_x})_{x\in
  X'})$ is a compatible triple, then for every $\ell\in G'$ and $x\in
X'$, with $i_\ell:h\mapsto \ell h\ell^{-1}$ the conjugation by $\ell$,
we have
$$
\nu_{G'_{\ell x}}=(i_\ell)_*\nu_{G'_x}\;.
$$
\elemm

\dem Let $x\in X'$, $\ell\in G'$, $H'=G'_x$ and $H''=G'_{\ell x}=\ell
H'\ell^{-1}$. Using the left invariance of $\nu_{X'}$ for the first
inequality and the bi-invariance of the Haar measure on $G'$ for the
last one, we have, for every $f:G'\ra\RR$ continuous with compact
support,
\begin{align*}
&\int_{g'\ell x\in X'}\int_{h''\in
  H''}f(g'h'')\;d(i_\ell)_*\nu_{G'_x}(h'')\;
d\nu_{X'}(g'\ell x)\\
=&\int_{g'\ell x\in X'}\int_{h''\in
  H''}f(\ell g'h'')\;d(i_\ell)_*\nu_{G'_x}(h'')\;
d\nu_{X'}(g'\ell x)\\ 
=&\int_{g'\ell x\in X'}\int_{h'\in
  H'}f(\ell g'\ell h'\ell^{-1})\;d\nu_{G'_x}(h')\;
d\nu_{X'}(g'\ell x)\\ 
=&\int_{g x\in X'}\int_{h'\in
  H'}f\circ i_\ell(gh')\;d\nu_{G'_x}(h')\;
d\nu_{X'}(g x)\\ =&\int_{G'} f\circ i_\ell\;d\nu_{G'}=\int_{G'} f\;d\nu_{G'}\;.
\end{align*}
The result then follows by uniqueness.
\cqfd

\medskip In order to deal with non simply connected groups, we
introduce a modified version of the Siegel weights.

Let ${\bf L}'$ be a reductive linear algebraic group defined and
anisotropic over $\QQ$, let $\pi:{\bf L}'\ra \GL({\bf V}')$ be a
rational linear representation of ${\bf L}'$ defined over $\QQ$, let
$v\in {\bf V}'(\QQ)$ be such that its orbit ${\bf X}'_v$ under ${\bf
  L}'$ is Zariski-closed in ${\bf V}'$, let $u\in {\bf X}'_v(\QQ)$ and
let ${\bf L}'_u$ be the stabiliser of $u$ in ${\bf L}'$. We define
the {\it modified Siegel weight} of $u$ as
\begin{equation}\label{eq:defimodifSiegelweight}
{w'}_{{\bf L}',\pi}(u)=
\frac{\vol\big( ({\bf L}'_u(\ZZ)\cap{\bf L}'(\RR)_0)\bs
({\bf L}'_u\cap{\bf L}'(\RR)_0)\big)}
{\vol\big( ({\bf L}'(\ZZ)\cap{\bf L}'(\RR)_0)\bs
{\bf L}'(\RR)_0\big)}\;,
\end{equation}
using Weil's convention for the normalisation of the measures on ${\bf
  L}'_u(\RR)$ (depending on the choice of a left Haar measure on ${\bf
  L}'(\RR)$ and of a ${\bf L}'(\RR)$-invariant measure on ${\bf
  X}'_v(\RR)$). Note that the denominator of the standard Siegel weight
$w_{{\bf L}',\pi}(u)$ is an integral multiple (depending only on ${\bf
  L}'$) of the denominator of the modified one, since $({\bf
  L}'(\ZZ)\cap{\bf L}'(\RR)_0)\bs {\bf L}'(\RR)_0$ is a connected
component of ${\bf L}'(\ZZ)\bs {\bf L}'(\RR)$.  But the ratio of the
numerator of the Siegel weight by the numerator of the modified one
may depend on $u$.

\medskip Let us now describe the framework of our main result. Let
${\bf G}$ be a connected reductive linear algebraic group defined over
$\QQ$. Let $\bf P$ be a (proper) parabolic subgroup of $\bf G$ defined
over $\QQ$ (see for instance \cite[\S III.1]{BorJi06}, \cite[\S
5.2]{Springer94})). Recall that a linear algebraic group defined over
$\QQ$ is {\it $\QQ$-anisotropic} if it contains no nontrivial
$\QQ$-split torus.

Recall that there exist a (nontrivial) maximal $\QQ$-split torus
${\bf S}$ in $\bf G$ (contained in $\bf P$ and unique modulo
conjugation by an element of ${\bf P}(\QQ)$), such that if
$\Phi^\CC=\Phi^\CC({\bf G},{\bf S})$ is the root system of $\bf G$
relative to ${\bf S}$ (seen contained in the set of characters of
${\bf S}$), if $\ggg^\CC_\beta$ is the root space of  $\beta\in
\Phi^\CC$, then there exist a unique set of simple roots
$\Delta=\Delta_{\bf P}$ in $\Phi^\CC$ and a unique proper subset
$I=I_{\bf P}$ of $\Delta$, such that, with $\Phi^\CC_+$ the set of
positive roots of $\Phi^\CC$ defined by $\Delta$ and $\Phi^\CC_I$ the
set of roots of $\Phi$ that are linear combinations of elements of
$I$, if ${\bf A}$ is the identity component of
$$
\bigcap_{\alpha\in I}\;\ker\alpha\;,
$$ 
which is a $\QQ$-split subtorus of $\bf S$, if $\bf U$ is the
connected algebraic subgroup of $\bf G$ defined over $\QQ$ whose Lie
algebra is
$$
\uuu^\CC=\bigoplus_{\beta\in \Phi^\CC_+-\Phi^\CC_I}\ggg^\CC_\beta\;,
$$
then $\bf P$ is the semi-direct product of its unipotent radical $\bf
U$ and of the centraliser of ${\bf A}$ in $\bf G$. Note that $\bf A$
is one-dimensional if $\bf P$ is a maximal (proper) parabolic subgroup
of $\bf G$ defined over $\QQ$ (that is, if $\Delta-I$ is a singleton).

Let $\ggg$ be the Lie algebra of ${\bf G}(\RR)$. Using the
multiplicative notation on the group of characters of $\bf S$, for
every $\alpha\in \Delta$, we define $m_\alpha=m_{\alpha,\bf P}\in\NN$
by
$$
\prod_{\beta\in \Phi^\CC_+-\Phi^\CC_I}
\;\beta^{\,\dim_\RR(\ggg^\CC_\beta\cap\ggg)}=
\prod_{\alpha\in \Delta}\;\alpha^{m_\alpha}\;.
$$
Let $(\alpha^\vee)_{\alpha\in \Delta-I}$ in ${{\bf A}
  (\RR)_0}^{\Delta- I}$ be such that $\log\beta (\alpha^\vee)$ is
equal to $1$ if $\alpha=\beta$ and to $0$ otherwise. Let
$\Lambda^\vee$ be the lattice in ${\bf A}(\RR)_0$ generated by
$\{\alpha^\vee\;:\;\alpha\in \Delta-I\}$. For every element
$T=(t_\alpha)_{\alpha\in\Delta-I}$ of $[0,+\infty[^{\Delta-I}$, let
$$
A_T=\{a\in{\bf A}(\RR)_0\;:\;\;
\forall\;\alpha\in\Delta- I,\;\;0\leq \log(\alpha(a))\leq t_\alpha\}\;.
$$

Recall that by the definition of {\it a relative Langlands
  decomposition } of the parabolic subgroup $\bf P$ defined over
$\QQ$, there exists a connected reductive algebraic subgroup $\bf M$
of $\bf P$ defined over $\QQ$ without nontrivial $\QQ$-characters such
that $\bf AM$ is the centraliser of $\bf A$ in $\bf G$. In particular,
$\bf AM$ is a Levi subgroup of $\bf P$ defined over $\QQ$, $\bf A$
centralises $\bf M$ and is the largest $\QQ$-split subtorus of the
centre of $\bf AM$, $\bf AM$ normalises $\bf U$, and
$$
{\bf P}={\bf AMU}\;.
$$

For every Lie group $G'$ endowed with a left Haar measure, for every
discrete subgroup $\Ga'$ of $G'$, we endow $\Ga'\bs G'$ with the
unique measure such that the canonical covering map $G'\ra \Ga'\bs G'$
locally preserves the measures.

In what follows, we will need a normalisation of the Haar measures,
which behaves appropriately when passing to some subgroups.  We will
start with a Riemannian metric on ${\bf G}(\RR)$, take the induced
Riemannian volumes on the real points of the various algebraic
subgroups of ${\bf G}$ defined over $\QQ$ that will appear, which will
give us the choices necessary for using Weil's normalisation to define
the Siegel weights.

\medskip
The main result of this paper is the following one.

\btheo \label{theo:main} Let $\bf G$ be a connected reductive linear
algebraic group defined over $\QQ$, without nontrivial
$\QQ$-characters. Let $G={\bf G} (\RR)_0$ and $\Ga={\bf G}(\ZZ)\cap
G$. Let $\bf P$ be a parabolic subgroup of $\bf G$ defined over $\QQ$,
and let ${\bf P} ={\bf AMU}$ be a relative Langlands decomposition of
$\bf P$.  Let $\rho:{\bf G}\ra \GL({\bf V})$ be a rational
representation of $\bf G$ defined over $\QQ$ such that there exists
$v_0\in{\bf V}(\QQ)$ whose stabiliser in $\bf G$ is ${\bf H}= {\bf
  MU}$.  Let $\bf L$ be a reductive algebraic subgroup of $\bf G$
defined and anisotropic over $\QQ$.

Assume that $\bf LP$ is Zariski-open in $\bf G$ and that for every
$a\in{\bf A}$, the orbit $\rho({\bf L} a)v_0$ is Zariski-closed in
$\bf V$. Endow ${\bf G} (\RR)$ with a left-invariant Riemannian
metric, for which the Lie algebras of ${\bf H}(\RR)$ and ${\bf
  A}(\RR)$ are orthogonal, and the orthogonal of the Lie algebra of
${\bf P}(\RR)$ is contained in the Lie algebra of ${\bf L}(\RR)$.

Then there exists $\delta>0$ such that, as $T=
(t_\alpha)_{\alpha\in\Delta-I}\in[0,+\infty[^{\Delta-I}$ and
$\min_{\alpha\in\Delta- I}t_\alpha$ tends to $+\infty$,
$$
\sum_{a\in A_T}\;\;\;
\sum_{[x]\in ({\bf L}(\RR)_0\cap \Ga)\bs(\rho({\bf
    L}(\RR)_0 a)v_0\cap \rho(\Ga)v_0)} \; 
{w'}_{{\bf L},\rho_{\mid {\bf L}}}(x)= 
$$
$$
\frac{\vol\big(({\bf H}\cap\Ga)\bs ({\bf H}\cap G)\big)
\vol(\Lambda^{\!\vee}\,\bs {\bf A}(\RR)_0)}
{\vol(\Ga\bs G)}\;
\Big(\prod_{\alpha\in \Delta- I}
\frac{e^{m_\alpha t_\alpha}}{m_\alpha}\Big)\big(1+
\operatorname{O}(e^{-\delta\min_{\alpha\in\Delta- I}t_\alpha})\big)\;.
$$
\etheo

\dem Let us start by fixing the notation that will be used throughout
the proof of Theorem \ref{theo:main}, and by making more explicit the
above-mentionned conventions about the various volumes that occur in
the asymptotic formula.

\medskip
Consider the connected real Lie group $G={\bf G}(\RR)_0$, its (closed)
Lie subgroups
$$
A={\bf A}(\RR)_0,\;H={\bf H}\cap G,\;L={\bf L}(\RR)_0,\;
M={\bf M}\cap G,\;P={\bf P}\cap G,\;U={\bf U}(\RR)\;.
$$
We have $H=MU$ and $P=AMU=MUA$, since $A$ and $U$ are connected. Note
that $L$ is also connected, but $H$ and $M$ are not necessarily
connected. We denote by 
$$
\aaa, \;\;\ggg,\;\; \hhh,\;\; \lll, \;\;\mmm,\;\;\ppp, \;\;\uuu
$$ 
the Lie algebras of the real Lie groups $A,G,H,L,M,P,U$
respectively, endowed with the restriction of the scalar product on
$\ggg$ defined by the Riemannian metric of $G$. 
Since $\bf L$ is $\QQ$-anisotropic, so is ${\bf L}\cap{\bf P}$. Since
the map ${\bf L}\cap{\bf P}\ra {\bf P}/{\bf H}\simeq{\bf A}$ is
  defined over $\QQ$ and $\bf A$ is a $\QQ$-split torus, this implies
  that the identity component of ${\bf L}\cap{\bf P}$ is contained in
  $\bf H$. In particular 
\begin{equation}\label{eq:lcappeglcaph}
\lll\cap\hhh=\lll\cap\ppp\;.
\end{equation}
Note that $\ggg=\lll+\ppp$ since ${\bf LP}$ is Zariski-open in $\bf
G$.  We have assumed that $\aaa$ is orthogonal to $\hhh$ and that the
orthogonal $\ppp^\perp$ of $\ppp$ is contained in $\lll$. In
particular, with $\qqq$ the orthogonal of $\lll\cap\hhh$ in $\hhh$,
we have the following orthogonal decompositions
\begin{equation}\label{eq:perpendicul}
\ggg=\ppp^\perp\stackrel{\perp}{\oplus}(\lll\cap\hhh)
\stackrel{\perp}{\oplus}\qqq\stackrel{\perp}{\oplus}\aaa,\;\;\;
\hhh=(\lll\cap\hhh)\stackrel{\perp}{\oplus}\qqq,\;\;\;
\lll=\ppp^\perp\stackrel{\perp}{\oplus}(\lll\cap\hhh),\;\;\;
\ppp=\hhh\stackrel{\perp}{\oplus}\aaa\;.
\end{equation}

The left-invariant Riemannian metric on $G$ induces a left Haar
measure $\omega_G$ on $G$, and a left-invariant Riemannian metric on
every Lie subgroup $G'$ of $G$, hence a left Haar measure
$$
\omega_{G'}
$$ 
on $G'$ (which is the counting measure if $G'$ is discrete). Note that
$A,G,H,L,M,U, L\cap H$ are unimodular: indeed $A,G,L,M$ are reductive
and $U$ is unipotent; furthermore, ${\bf L}\cap {\bf H}$ is the
stabiliser of $v_0$ in ${\bf L}$, the orbit of $v_0$ under $\bf L$ is
affine and hence ${\bf L}\cap {\bf H}$ is reductive by
\cite[Theo.~3.5]{BorelHarishChandra62}. But $P$ is not unimodular.

The map $A\times M\times U\ra P$ defined by $(a,m,u)\mapsto amu$ is a
smooth diffeomorphism (see for instance \cite[page 273]{BorJi06}).  We
will denote by $d\omega_Ad\omega_H$ the measure on $P$ which is the
push-forward of the product measure by the diffeomorphism $(a,h)
\mapsto ah$.  Since $A$ normalises $H$, the measure $d\omega_A
d\omega_H$ is left-invariant by $P$, so that $d\omega_P(ah)$ and
$d\omega_A(a) d\omega_H(h)$ are proportional. Since these measures are
induced by Riemannian metrics, and since $\aaa$ and $\hhh$ are
orthogonal, we hence have
$$
d\omega_P(ah)=d\omega_A(a)d\omega_H(h)\;.
$$
Since $A$ normalises $U$, the group $A$ acts on the Lie algebra $\UUU$
of $U$ by the adjoint representation. The roots of this linear
representation of $A$ are exactly the restrictions to $A$ of the
elements $\beta$ in $\Phi^\CC_+ - \Phi^\CC_I$, with root spaces
$\ggg^\CC_\beta \cap\ggg$ and a set of simple roots is the set of
restrictions of the elements of $\Delta-I$ to $A$ (see for instance
\cite[Rem.~III.1.14]{BorJi06}).  Since $A$ is connected, these roots
have value in $]0,+\infty[\,$. The map $A\ra \RR^{\Delta-I}$ defined
by $a\mapsto (\log(\alpha(a)))_{\alpha\in \Delta-I}$ is hence a smooth
diffeomorphism.  We will denote by $\prod_{\alpha\in\Delta-I}
dt_\alpha$ the measure on $A$ which is the push-forward of the product
Lebesgue measure by the inverse of this diffeomorphism.  By
invariance, there exists a constant $c_A>0$ such that
$$
d\omega_A=c_A\prod_{\alpha\in\Delta-I}dt_\alpha\;.
$$
By the definition of $\Lambda^\vee$, we have $c_A=\Vol(\Lambda^{\vee}\bs
A)$.

Let $\Ga={\bf G}(\ZZ)\cap G$, which is a discrete subgroup of $G$
acting isometrically for the Riemannian metric of $G$ by left
translations. Let $Y_G=\Ga\bs G$ and let $\pi:G\ra Y_G=\Ga\bs G$ be
the canonical projection, which is equivariant under the right actions
of $G$. Then $Y_G$ is a connected Riemannian manifold (for the unique
Riemannian metric such that $\pi$ is a local isometry) endowed with
the transitive right action of $G$ by translations on the right. To
simplify the notation, for every Lie subgroup $G'$ of $G$, define
$$
Y_{G'}=\pi(G')\;,
$$
which is a injectively immersed submanifold in $Y_G$, endowed with the
Riemannian metric induced by $Y_G$, and identified with $(G'\cap \Ga)\bs
G'$ by the map induced by the inclusion of $G'$ in $G$.  Note that
$Y_L$ and $Y_U$ are connected, but $Y_H$ and $Y_M$ are not necessarily
connected. For every Lie subgroup $G'$ of $G$, let
$$
\mu_{G'}
$$ 
be the Riemannian measure on $Y_{G'}$, which locally is the
push-forward of the left Haar measure $\omega_{G'}$.

Since $\bf G$ and the identity component of $\bf MU$ have no
nontrivial $\QQ$-character, the Riemannian manifolds $Y_G$ and $Y_H$
have finite volume (see \cite[Theo.~9.4]{BorelHarishChandra62}) and
$Y_H$ is closed in $Y_G$ (see for instance
\cite[Theo.~1.13]{Raghunathan72}). Since $\bf L$ is reductive and
$\QQ$-anisotropic, the submanifold $Y_L$ is compact (see
\cite[Theo.~11.6]{BorelHarishChandra62}). Since $\bf U$ is unipotent,
the submanifold $Y_U$ is compact (see for instance
\cite[\S~6.10]{BorelHarishChandra62}).

For every Lie subgroup $G'$ of $G$ such that $Y_{G'}$ has finite
measure (that is, such that $\Ga\cap G'$ is a lattice in $G'$), we
denote by
$$
\overline\mu_{G'}=\frac{\mu_{G'}}{\|\mu_{G'}\|}
$$
the finite measure $\mu_{G'}$ normalised to be a probability
measure. In particular, $\overline\mu_G$, $\overline\mu_H$,
$\overline\mu_L$, $\overline\mu_U$ are well defined.

For every $T=(t_\alpha)_{\alpha\in\Delta-I}$ and $T'=
(t'_\alpha)_{\alpha\in\Delta-I}$ in $[0,+\infty[^{\Delta-I}$, let
$$
A_{[T,T']}=\{a\in A\;:\;\; \forall\;\alpha\in\Delta- I,
\;\;t_\alpha\leq \log(\alpha(a))\leq t'_\alpha\}\;.
$$
and $P_{[T,T']}=UM{A_{[T,T']}}^{-1}=H{A_{[T,T']}}^{-1}$. Define
$Y_{P_{[T,T']}}= \pi(P_{[T,T']})$, which is a submanifold with
boundary of $Y_G$, invariant under the right action of $H$, since $A$
normalises $H$.  To shorten the notation, we define
$$
A_T=A_{[0,T]}\;,\;\;\;P_T=P_{[0,T]}=H{A_T}^{-1}\;\;\;{\rm and}\;\;\;\;
Y_{P_T}=Y_{P_{[0,T]}}=\pi(P_T)=Y_HA_T^{-1}\;,
$$
as well as $\min T=\min_{\alpha\in \Delta-I}t_\alpha\geq 0$, which measures
the complexity of $T$ and will converge to $+\infty$. We will need to
estimate the volume of $\pi(P_{T})$ for $\mu_P$.

\blemm\label{lem:croissancePT} For every
$T=(t_\alpha)_{\alpha\in\Delta-I}$ in $[0,+\infty[^{\Delta-I}$, we
have
$$
\mu_P(Y_{P_{T}})=\vol(\Lambda^{\!\vee}\,\bs {\bf A}(\RR)_0)\;\;\|\mu_H\|
\prod_{\alpha\in \Delta-I}\frac{e^{m_\alpha t_\alpha}}{m_\alpha}\;.
$$
\elemm

\dem Denote by $du_\beta$ the Lebesgue measure on the Euclidean space
$\ggg_\beta=\ggg_\beta^\CC\cap\ggg$. For any order on $\Phi^\CC_+ -
\Phi^\CC_I$, the map from $\prod_{\beta\in \Phi^\CC_+ - \Phi^\CC_I}
\;\ggg_\beta$ to $U$ defined by $(u_\beta)_{\beta\in \Phi^\CC_+ -
  \Phi^\CC_I} \mapsto \prod_{\beta\in \Phi^\CC_+ - \Phi^\CC_I}\exp
u_\beta$ is a smooth diffeomorphism, and there exists $c_U>0$ such
that $d\omega_U$ is the push-forward by this diffeomorphism of the
measure $c_U\prod_{\beta\in \Phi^\CC_+ - \Phi^\CC_I}du_\beta$.

For every $a\in A$, if $i_a:g\mapsto aga^{-1}$ is the conjugation by
$a$, then for every $u_\beta\in\ggg_\beta$, we have $i_a(\exp u_\beta)
= \exp ((\Ad a)(u_\beta))=\exp (\beta(a)u_\beta)$. Hence
$$
(i_a^{-1})_*(\omega_U)=\prod_{\beta\in \Phi^\CC_+ -
  \Phi^\CC_I}\;\beta(a)^{\dim \ggg_\beta}\omega_U
=\prod_{\alpha\in \Delta -I}\;\alpha(a)^{m_\alpha}\omega_U
$$
by the definition of $(m_\alpha)_{\alpha\in \Delta}$ and since the
elements of $I$ are trivial on $A$. Since $A$ commutes with $M$, we
hence have $(i_a^{-1})_*(\omega_H)=\prod_{\alpha\in \Delta -I}
\;\alpha(a)^{m_\alpha}\omega_H$.

We have, since $A$ is unimodular,
$$
d\omega_P(ha^{-1})=
d\omega_P(a^{-1}aha^{-1})=d\omega_A(a^{-1})d\omega_H(aha^{-1})
=d\omega_A(a)d((i_a^{-1})_*\omega_H)(h)\;.
$$
Since $\Ga\cap P=\Ga\cap H$ (see for instance the lines following
Proposition III.2.21 in \cite[page 285]{BorJi06}) and $A\cap H=\{e\}$,
we have $\pi(Ha)\neq \pi(Ha')$ if $a\neq a'$. Hence
\begin{align}
  \mu_P(Y_{P_{T}})&=\int_{y\in Y_H}\int_{a\in
    A_{T}}d\mu_P(ya^{-1})= \int_{A_{T}}\prod_{\alpha\in \Delta
    -I} \;\alpha(a)^{m_\alpha}d\omega_A(a)\int_{Y_H}d\mu_H
\label{eq:dandemo}\\ &
  =\|\mu_H\|\,c_A\prod_{\alpha\in \Delta-I}\int_{0
  }^{t_\alpha}e^{m_\alpha s}\;ds\:.\nonumber
\end{align}
Since $m_\alpha>0$, the result follows. 
 \cqfd

\medskip
To simplify the notation, we write $\rho(g)x=gx$ for every $g \in {\bf
  G}$ and $x\in{\bf V}$, we define $v_a=av_0$ for every $a\in A$, and
we denote by $L_x={\bf G}_x\cap L$ the stabiliser of $x$ in $L$
for every $x\in {\bf V}(\RR)$. 

Since we have a left Haar measure $\omega_L$ on $L$ and $\omega_{L\cap
  H}$ on $L\cap H$, Weil's normalisation gives a $L$-invariant measure
on the homogeneous space $L/(L\cap H)$, and hence a left Haar measure
on the stabilisers $\ell(L\cap H) \ell^{-1}$ for every $\ell$ in $L$,
as explained above.  As announced, the modified Siegel weights
${w'}_{{\bf L}, \rho_{\mid {\bf L}}}(\cdot)$ are defined using this
Weil's normalisation, as follows.

For every $\ell\in L$ and $a\in A$, if $x=\ell av_0$,  since $H=MU$ is
normalised by $A$ and is the stabiliser of $v_0$ in $G$, we have
\begin{equation}\label{eq:calclsubx}
L_x=L\cap\operatorname{Stab}_Gx= L\cap (\ell H\ell^{-1})=
\ell (L\cap H)\ell^{-1}\;.
\end{equation}
Note that 
\begin{equation}\label{eq:interintersimpl}
{\bf L}(\RR)_0\cap {\bf L}(\ZZ)= {\bf L}(\RR)_0\cap {\bf G} (\ZZ)=
{\bf L}(\RR)_0\cap{\bf G}(\RR)_0\cap {\bf G} (\ZZ) =L\cap\Ga\;,
\end{equation}
and similarly ${\bf L}(\RR)_0\cap {\bf L}_x(\ZZ)=L_x\cap\Ga$ for every
$x\in Lv_a\cap\Ga v_0$.  Hence the denominator of the modified Siegel
weight ${w'}_{{\bf L}, \rho_{\mid {\bf L}}}(x)$ is equal to $\vol
\big((L\cap\Ga) \bs L\big) =\vol(Y_L)$, using the measure $\mu_L$ on
$Y_L$ induced by the Haar measure $\omega_L$ on $L$. Its numerator
$\vol\big((L_x\cap\Ga)\bs L_x\big)$ is defined by using the measure on
$(L_x\cap\Ga)\bs L_x$ induced by the left Haar measure on $L_x=\ell
(L\cap H)\ell^{-1}$ given by Weil's normalisation.

\medskip
We now proceed to the proof of Theorem
\ref{theo:main}.

\bigskip 
\noindent{\bf Step 1. } 
The first step of the proof is the following group theoretic lemma,
which relates the counting function of modified Siegel weights to the
counting function of volumes of orbits of $L\cap H$.  We denote with
square brackets the (left or right) appropriate orbit of an element.

\blemm \label{lem:grouptheoretic}
For every $a\in A$, there exists a bijection between finite
subsets
$$
\Theta_a\,:\;(L\cap \Ga)\bs (Lv_a\cap\Ga v_0)\;\longrightarrow\;
(Y_L\cap Y_H a^{-1})/(L\cap H)
$$
such that for every $x\in Lv_a\cap\Ga v_0$, if $[y]=\Theta_a([x])$,
then
\begin{equation}\label{eq:equalvolsiegel}
\vol\big((L_x\cap \Ga)\bs L_x\big)=\vol\big(y(L\cap H)\big)\;.
\end{equation}
In particular, for every $T\in[0,+\infty[^{\Delta-I}$, we have
$$
\sum_{a\in A_{T}}\;\;
\sum_{[x]\in (L\cap \Ga)\bs(Lv_a\cap \Ga v_0)} 
\;  {w'}_{{\bf L},\rho_{\mid {\bf L}}}(x)= \sum_{a\in A_{T}}\;\;
\sum_{[y]\in (Y_L\cap Y_H a^{-1})/(L\cap H)} \;
\frac{\vol\big(y(L\cap H)\big)}{\vol(Y_L)}\;. 
$$
\elemm

\dem Fix $a\in A$. First note that the groups $L\cap \Ga$ and $L\cap
H$ do preserve the subsets $Lv_a\cap\Ga v_0$ of $V(\RR)$ and $Y_L\cap
Y_H a^{-1}$ of $Y$ respectively for their left and right action, since
$A$ normalises $H$. The finiteness of the set $(L\cap \Ga)\bs
(Lv_a\cap\Ga v_0)$ follows from Borel-Harish-Chandra's finiteness
theorem as in the introduction. Also recall that $H=MU$ is the
stabiliser of $v_0$ in $G$.

Define 
$$
\Theta_a:[\ell v_a]\mapsto [\pi(\ell)]\;.
$$
Let us prove that this map is well defined and bijective. Let
$\ell,\ell'\in L$.

We have $\ell v_a\in\Ga v_0$ if and only if there exists $\ga\in\Ga$
such that $\ell av_0=\ga v_0$, that is, if and only if there exist
$\ga\in\Ga$ and $h\in H$ such that $\ell =\ga ha^{-1}$, that is, if and
only if $\pi(\ell)\in Y_L\cap Y_Ha^{-1}$. This proves that $\Theta_a$
has values in $(Y_L\cap Y_H a^{-1})/(L\cap H)$ and is surjective.

Let us prove that $\Theta_a$ does not depend on the choice of
representatives and is injective. We have $[\ell' v_a]=[\ell v_a]$ if
and only if there exists $\ga\in L\cap \Ga$ such that $\ell' a v_0=\ga
\ell a v_0$, hence if and only if there exist $\ga\in L\cap \Ga$ and
$h\in H$ such that $\ell' a =\ga \ell a h$. Note that this equation
implies that $aha^{-1}\in L$ if and only if $\ga\in L$. Since $A$
normalises $H$, we hence have $[\ell' v_a]=[\ell v_a]$ if and only if
$\ell'\in\Ga\ell (L\cap H)$, that is if and only if
$[\pi(\ell)]=[\pi(\ell')]$.

\medskip To prove the second assertion, let $\ell\in L$ be such that
$x=\ell v_a\in\Ga v_0$ and let $y= \pi(\ell)$, so that
$\Theta_a([x])=[y]$. The orbit of $y$ under $L\cap H$ in $Y_G$ is the
image by the locally isometric map $\pi$ of the Riemannian submanifold
$\ell (L\cap H)$ of $G$. The left translation by $\ell^{-1}$ is an
isometry (hence is volume preserving) from $\ell (L\cap H)$ to $(L\cap
H)$. By Lemma \ref{lem:weilnormal}, the map $g\mapsto \ell g\ell^{-1}$
from $L\cap H$ to $\ell(L\cap H)\ell^{-1}$, which is equal to $L_x $
by Equation \eqref{eq:calclsubx}, is measure preserving.  Therefore
the map $\varphi:[z]\mapsto [z\ell]$ from $(L_x\cap \Ga)\bs L_x$ to
$y(L\cap H)$ is a measure preserving bijection. This proves the volume
equality of Equation \eqref{eq:equalvolsiegel}.

The last claim follows from the other ones, since the numerator of the
modified Siegel weight $w'_{{\bf L},\rho_{\mid \bf L}}(x)$ is
$\vol\big((L_x\cap \Ga)\bs L_x\big)$.  \cqfd

\bigskip 
\noindent{\bf Step 2. } 
The second step of the proof is an equidistribution result, in the
spirit of \cite{KleMar96}, saying that the piece $Y_{P_T}$ of orbit of
$P$ equidistributes in $Y_G$ as $\min T\ra +\infty$.

For every smooth Riemannian manifold $Z$ and $\ellk\in\NN$, we denote
by $\C_c^\ellk(Z)$ the normed vector space of $\operatorname{C}^\ellk$
maps with compact support on $Z$, with norm $\|\cdot\|_\ellk$.

\bprop\label{prop:step2} There exist $\ellk\in\mathbb{N}$ and
$\kappa>0$ such that for every $f\in \C_c^\ellk(Y_G)$ and
$T=(t_\alpha)_{\alpha\in \Delta-I}\in [0,+\infty[^{\Delta-I}$, we
have, as $\min T$ tends to $+\infty$,
$$
\frac{1}{{\mu_P(Y_{P_T})}}\int_{Y_{P_T}} f\,d\mu_P =\int_{Y_G} f\,
d\overline\mu_G + \operatorname{O}\big(e^{-\kappa\min T} \;\|f\|_\ellk
\big)\;.
$$
\eprop

To prove the proposition, we will use the disintegration formula already
seen in the proof of Lemma \ref{lem:croissancePT}
\begin{align*}
\int_{Y_{P_T}} f\,d\mu_P
=\int_{A_T}\Big(\int_{Y_H} f(ya^{-1})\, d\mu_H(y)\Big)\,
\Big(\prod_{\alpha\in \Delta-I} \alpha(a)^{m_\alpha}\Big)d\omega_A(a)\;.
\end{align*}
This formula indicates that the proposition would follow from 
(an averaging of) the equidistribution of the translates $Y_Ha^{-1}$,
which is established in Proposition \ref{th:equidist_H} below.
To state this proposition, we need to introduce additional notation.

The linear algebraic group ${\bf G}$ decomposes as an almost direct
product
$$
{\bf G}={\bf Z}({\bf G}){\bf G}_1\cdots {\bf G}_s
$$
where ${\bf Z}({\bf G})$ is the centre of ${\bf G}$, and ${\bf
  G}_1,\ldots, {\bf G}_s$ are $\mathbb{Q}$-simple connected algebraic
subgroups of ${\bf G}$.  The maximal $\QQ$-split torus ${\bf S}$
decomposes as an almost direct product
$$
{\bf S}={\bf S}_1\cdots {\bf S}_s
$$
where ${\bf S}_i$ is a maximal $\QQ$-split torus in ${\bf G}_i$. We
also get an almost direct product decomposition
\begin{equation}\label{eq:decomp_G}
G=Z(G)G_1\cdots G_s,
\end{equation}
where $Z(G)$ is the centre of $G$ (which is equal to ${\bf Z}({\bf
  G})_0$ since $\bf G$ is connected and $G$ is Zariski-dense in $\bf
G$) and $G_i={\bf G}_i(\mathbb{R})_0$ for $1\leq i\leq s$.  Since $\bf
G$ has no nontrivial $\QQ$-character, and since $\bf M$ is the
centraliser of $\bf A$, this gives corresponding almost direct product
decompositions of the Lie groups $A=A_1\dots A_s$ (this one being a
direct product), $U=U_1\dots U_s$, $M=Z(G)M_1\dots M_s$,
$H=Z(G)H_1\dots H_s$.  The set of simple roots $\Delta$ decomposes as
a disjoint union
$$
\Delta=\Delta_1\sqcup\cdots \sqcup\Delta_s
$$
where $\Delta_i$ is a set of simple roots of ${\bf G}_i$ relatively to
${\bf S}_i$, and the positive (closed) Weyl chamber $A^+$ in $A$
associated to $\Delta$ decomposes as
$$
A^+=A^+_1\cdots  A^+_s,
$$
where $A^+_i$ is the positive (closed) Weyl chamber in $A\cap
G_i$ associated to $\Delta_i$.

For $1\leq i\leq s$ and $a\in A_i$, we define
\begin{equation}\label{eq:E_i}
E_i(a)=\exp\big(-\big(\max_{\alpha\in\Delta_i - I}
      \log\alpha(a)\big)\big)>0
\end{equation}
if $\Delta_i - I\ne \emptyset$, and $E_i(a)=0$, otherwise. For every
$\kappa>0$, we also define
$$
E^\kappa(a)=\sum_{i=1}^s E_i(a_i)^\kappa
$$ 
for every $a\in A^+$ with $a_1\in A_1^+, \cdots, a_s\in A_s^+$ and
$a=a_1\cdots a_s$.

\bprop\label{th:equidist_H} There exist $\ellk\in\mathbb{N}$ and
$\kappa>0$ such that for all $f\in \C_c^\ellk(Y_G)$ and $a\in A^+$,
$$
\int_{Y_H} f(ya^{-1})\, d\overline\mu_H(y)
=\int_{Y_G} f\, d\overline\mu_G+
\operatorname{O}\big(E^\kappa(a)\;\|f\|_\ellk\big)\;.
$$
\eprop

Given a  Lie subgroup $D$ of $G$ such that $\Gamma\cap D$ is a
lattice in $D$, we denote by $\overline\nu_D$ the normalised right
invariant measure on $(\Gamma\cap D)\backslash D$. Recall that
$Y_D=\pi(D)$ is a closed submanifold of $Y_G$, and that $\mu_D$ is the
invariant measure on $Y_D$ induced by the Riemannian metric, with
normalised measure $\overline\mu_D$.

We identify $(\Gamma\cap H)\backslash H$ with $Y_H$ using the (well
defined) map $h\mapsto \Gamma h$ (denoting again by $h\in H$ a
representative of a coset $h\in (\Gamma\cap H)\backslash H$).  Since
the groups $Z(G),H_1,\cdots, H_s$ commute, we also have the map
$$
(\Gamma\cap Z(G))\backslash
Z(G)\times (\Gamma\cap H_1)\backslash H_1\cdots \times (\Gamma\cap
H_s)\backslash H_s \to Y_H
$$ 
well defined by $(h_0,h_1,\ldots,h_s)\mapsto \Gamma h_0h_1\cdots h_s$
(using conventions similar to the above one for coset
representatives). Then the normalised invariant measures
$\overline\mu_H$, $\overline\mu_{H_1}, \dots, \overline\mu_{H_s}$
satisfy, for all $f\in \C_c(Y_G)$,
\begin{align}
&\int_{Y_H} f(y)\, d\overline\mu_H(y)=
\int_{(\Gamma\cap H)\backslash H} f(\Gamma h)
d\overline\nu_H(h)\nonumber\\
=&
\int_{(\Gamma\cap Z(G))\backslash Z(G)\times \cdots
\times (\Gamma\cap H_s)\backslash H_s} f(\Gamma h_0h_1\cdots h_s)
d\overline\nu_{Z(G)}(h_0)\cdots d\overline\nu_{H_s}(h_s).\label{eq:Y_H}
\end{align}
We will prove Proposition \ref{th:equidist_H} by using an inductive
argument on the number of factors.  We start by analysing the
distribution of $Y_{U_i}a^{-1}$ in Lemma \ref{lem:step_U} and then the
distribution of $Y_{H_i}a^{-1}$ in Lemma \ref{lem:step_H}.

Let $D$ be a product of almost direct factors of $G$ in the
decomposition \eqref{eq:decomp_G}.  For every $f\in\C_c^0(Y_G)$, we
define a map $\mathcal{P}_D f:Y_G\ra\CC$ by
$$
(\mathcal{P}_D f)(\Gamma g)=\int_{(\Gamma\cap D)\backslash D} f(\Gamma d g)\;
d\overline\nu_D(d)
$$
which does not depend on the choice of the representative of $\Ga g$,
by the right invariance of $\overline\nu_D$ under $D$. Note that
$\mathcal{P}_D f$ is continuous and invariant under the right action
of $D$.

\blemm\label{lem:step_U}
There exist $\ellk\in\mathbb{N}$ and $\kappa_1>0$ such that 
for every $i\in\{1,\dots,s\}$ with $A_i\ne \{1\}$, 
for every $f\in \C_c^\ellk(Y_G)$ and $a\in A_i^+$,
\begin{align*}
\int_{Y_{U_i}} f(ya^{-1})\, d\overline\mu_{U_i}(y)
=(\mathcal{P}_{G_i} f)(\Gamma e)+
\operatorname{O}\big(E_i(a)^{\kappa_1}\|f|_{Y_{G_i}}\|_{\ellk}\big)\;.
\end{align*}
\elemm

\dem For $1\leq i\leq s$, we consider the unitary representation of
the group $G_i$ on the orthogonal complement of the space of
$G_i$-invariant (hence constant on $Y_{G_i}$) functions in the Hilbert
space $\LL^2(Y_{G_i}, \overline{\mu}_{G_i})$, whose scalar product we
denote by $\langle\cdot,\cdot\rangle_{Y_{G_i}}$ (using the normalised
measure $\overline{\mu}_{G_i}$).  We note that for every $f\in
\C_c(Y_G)$, the function $f|_{Y_{G_i}} -(\mathcal{P}_{G_i}f)(\Gamma
e)$ belongs to this space.

We say that a unitary representation of a connected real semisimple
Lie group $G'$ has the {\it strong spectral gap property} if the
restriction to every noncompact simple factor of $G'$ is isolated from
the trivial representation for the Fell topology (see for instance
\cite{Cowling79}, \cite[Appendix]{BekHarVal08},
\cite[Appendix]{KleMar99} for equivalent definitions and examples, and
compare for instance with \cite{Nevo05,KelSar09} for variations on the
terminology).  We claim that the above unitary representation of $G_i$
has the strong spectral gap property.  Indeed, if ${\bf G}_i$ is
simply connected and $\Gamma_i$ is a congruence subgroup in $G_i$,
then the strong spectral gap property on $\Gamma_i\backslash G_i$ is a
direct consequence of the property $\tau$ proved in \cite{Clozel03},
see Theorem 3.1 therein.  By \cite[Lemma~3.1]{KleMar99}, this also
implies, when ${\bf G}_i$ is simply connected, the strong spectral
property on $\Gamma_i\backslash G_i$ for subgroups $\Gamma_i$ that are
commensurable with congruence subgroups, and, in particular, for
arithmetic subgroups of $G_i$.  Now let $p_i:\widetilde {\bf G}_i\to
{\bf G}_i$ be a simply connected cover of ${\bf G}_i$, and let
$\widetilde G_i=\widetilde{\bf G}_i(\mathbb{R})$.  Then $Y_{G_i}\simeq
p_i^{-1}(\Gamma\cap G_i)\backslash\widetilde G_i$, and the strong
spectral gap property for $\LL^2(Y_{G_i}, \overline{\mu}_{G_i})$
follows from the above arguments.
 
Applying \cite[Theorem~2.4.3]{KleMar96}, we deduce that there exist
$\ellk\in\NN$ and $\kappa_1',C>0$ such that for every $i\in\{1,\dots,s\}$
such that $A_i\neq \{e\}$, for every $\phi\in \C_c^\ellk(Y_{G_i})$ and
$a\in A_i^+$,
\begin{equation}\label{eq:mixing}
\big\langle(f|_{Y_{G_i}} -(\mathcal{P}_{G_i}f)(\Gamma e))\circ a^{-1}, 
\phi\,\big\rangle_{Y_{G_i}}\leq C\,
E_i(a)^{\kappa_1'} \;\|f|_{Y_{G_i}}\|_{\ellk}\;\|\phi\|_{\ellk}\;,
\end{equation}
where $E_i(a)$ is defined in Equation \eqref{eq:E_i}.

Let $P^-_i$ denote the parabolic subgroup in $G_i$ opposite to $U_i$.
The product map $U_i\times P_i^-\to G_i$ is a diffeomorphism between
neighbourhoods of the identities. Since $Y_{U_i}=\pi(U_i)$ is compact,
if $\epsilon>0$ is small enough, there exists an open
$\epsilon$-neighbourhood $\Omega_\epsilon$ of the identity in $P^-_i$
such that the product map $Y_{U_i}\times \Omega_\epsilon \to Y_{G_i}$
is a diffeomorphism onto its image $Y_{U_i} \Omega_\epsilon$. We have
(see also Lemma \ref{lem:desintegr})
$$
\forall\;y\in Y_{U_i}, \;\forall\;p\in \Omega_\epsilon, 
\quad \quad d\overline\mu_{G_i}(yp)=d\overline\mu_{U_i}(y)d\omega(p),
$$
for a suitably normalised smooth measure $\omega$ on
$\Omega_\epsilon$.  There exists $\sigma>0$ (depending on $\ellk$) such
that for every $\epsilon>0$ small enough, there exists a nonnegative
function $\psi_\epsilon\in \C^\ellk_c(\Omega_\epsilon)$ satisfying
$$
\int_{\Omega_\epsilon} \psi_\epsilon \, d\omega=1\quad\hbox{and}\quad
\|\psi_\epsilon\|_\ellk=\operatorname{O}(\epsilon^{-\sigma})\;.
$$
Define a $\operatorname{C}^\ellk$ function $\phi_\epsilon:Y_{G_i}
\ra[0,+\infty[$ supported on $Y_{U_i}\Omega_\epsilon$ by
$$
\forall\;y\in Y_{U_i}, \;\forall\;p\in \Omega_\epsilon, 
\quad\quad \phi_\epsilon(yp)=\psi_\epsilon(p)\;.
$$
Then
$$
\int_{Y_{G_i}} \phi_\vre \, d\overline\mu_{G_i}=1\quad\hbox{and}\quad
\|\phi_\epsilon\|_\ellk=\operatorname{O}(\epsilon^{-\sigma})\;.
$$
Since for all $a\in A_i^+$ and $p\in \Omega_\epsilon$,
$$
d(apa^{-1},e)=\operatorname{O}(\epsilon),
$$
we obtain
\begin{align*}
  \big\langle f|_{Y_{G_i}}\circ a^{-1},
  \phi_\epsilon\big\rangle_{Y_{G_i}} &=\int_{Y_{U_i}\times
    \Omega_\epsilon} f(ypa^{-1})
  \psi_\epsilon(p)\, d\overline\mu_{U_i}(y)d\omega(p)\\
  &=\int_{Y_{U_i}} f(ya^{-1})\,
  d\overline\mu_{U_i}(y)+\operatorname{O}\left(\epsilon
    \|f|_{Y_{G_i}}\|_1\right).
\end{align*}
Since $\mathcal{P}_{G_i}f$ is $G_i$-invariant,
\begin{align*}
  \big\langle\mathcal{P}_{G_i}f, \phi_\epsilon\big\rangle_{Y_{G_i}} =
  (\mathcal{P}_{G_i}f)(\Gamma e)\Big(\int_{Y_{G_i}} \phi_\vre \,
    d\overline\mu_{G_i}\Big) = (\mathcal{P}_{G_i}f)(\Gamma e).
\end{align*}
Combining these estimates with \eqref{eq:mixing} (we may assume
that $\ellk\geq 1$), we conclude that
\begin{align*}
  \int_{Y_{U_i}} f(ya^{-1})\, d\overline\mu_{U_i}(y)
  =(\mathcal{P}_{G_i} f)(\Gamma e) +\operatorname{O}\big((\epsilon
    +E_i(a)^{\kappa_1'} \epsilon^{-\sigma})\;
    \|f|_{Y_{G_i}}\|_\ellk\big).
\end{align*}
Finally, taking $\epsilon=E_i(a)^{\kappa_1'/(1+\sigma)}$ which is
small if $a$ lies outside a compact subset of $A^+_i$, we deduce that
\begin{align*}
  \int_{Y_{U_i}} f(ya^{-1})\, d\overline\mu_{U_i}(y)
  =(\mathcal{P}_{G_i} f)(\Gamma e)
  +\operatorname{O}\left(E_i(a)^{\kappa_1'/(1+\sigma)}
    \|f|_{Y_{G_i}}\|_\ellk\right),
\end{align*}
as required.
\cqfd

\blemm\label{lem:step_H} There exist $\ellk\in\mathbb{N}$ and
$\kappa_2>0$ such that for every $f\in \C_c^\ellk(Y_G)$ and $a\in
A_i^+$, for every $i\in\{1,\dots,s\}$, we have
\begin{align*}
\int_{Y_{H_i}} f(ya^{-1})\, d\overline\mu_{H_i}(y)
=(\mathcal{P}_{G_i} f)(\Gamma e)+
\operatorname{O}\big(E_i(a)^{\kappa_2}\|f|_{Y_{G_i}}\|_{\ellk}\big).
\end{align*}
\elemm

\dem We first observe that if $A_i=\{e\}$, then $H_i=G_i$, and the
claim of the lemma is obvious.  Now we assume that $A_i\ne \{e\}$ in
which case Lemma \ref{lem:step_U} applies.

Let $N_i=(\Gamma\cap M_i)\backslash M_i$.  The space $Y_{H_i}=
\pi(U_iM_i)$ is a bundle over $N_i$ with fibres isomorphic to
$Y_{U_i}$, and the invariant measure $\overline\mu_{H_i}$ on $Y_{H_i}$
decomposes with respect to this structure.  Explicitly, for every
$m\in M_i$, the integrals $\int_{Y_{U_i}} f(ym)\;
d\overline\mu_{U_i}(y)$ for all $f\in \C_c(Y_G)$ define a
$U_i$-invariant probability measure on $Y_{U_i}m$, which depends only
on the coset $n=[m]$ of $m$ in $N_i=(\Gamma\cap M_i)\backslash M_i$,
and the $H_i$-invariant probability measure on $Y_{H_i}$ is given by $
\int_{N_i} \big(\int_{Y_{U_i}} f(ym)\, d\overline\mu_{U_i}(y)
\big)d\overline\nu_{M_i}([m])$ for all $f\in \C_c(Y_G)$.  Hence,
denoting again by $n$ any representative of a coset $n$ in $N_i$,
since $A$ centralises $M$,
\begin{align*}
\int_{Y_{H_i}} f(ya^{-1})\, d\overline\mu_{H_i}(y)
&=\int_{N_i}\int_{Y_{U_i}} f(yna^{-1})\; d\overline\mu_{U_i}(y)\;
d\overline\nu_{M_i}(n)\\
&=\int_{N_i}\int_{Y_{U_i}} f(ya^{-1}n)\, d\overline\mu_{U_i}(y)
\;d\overline\nu_{M_i}(n)\;.
\end{align*}

For $m\in M_i$ and $f\in\C_c(Y_G)$, we consider the function $f_m:
Y_G\ra\CC$ defined by $y\mapsto f(ym)$.  We note that there exist
$c_1,C'>0$ such that for every $f\in\C_c^\ellk(Y_G)$, we have
$$
\|f_m|_{Y_{G_i}}\|_{\ellk}\leq C' e^{c_1 d(e,m)} \|f|_{Y_{G_i}}\|_{\ellk}\;.
$$
Hence, by Lemma \ref{lem:step_U}, for every $f\in \C_c^\ellk(Y_G)$ and
$m\in M_i$, since $\mathcal{P}_{G_i} f_m=\mathcal{P}_{G_i} f$ by
invariance under $G_i$,
\begin{align}\label{eq:km2}
  \int_{Y_{U_i}} f(ya^{-1}m)\,d\overline\mu_{U_i}(y)=
  (\mathcal{P}_{G_i} f)(\Gamma
  e)+\operatorname{O}\big(E_i(a)^{\kappa_1} e^{c_1\,
    d(e,m)}\|f|_{Y_{G_i}}\|_\ellk\big)\;.
\end{align}
We fix $n_0\in N_i$ and for $R>0$, we set
$$
(N_i)_R=\{n\in N_i:\, d(n_0,n)\le R\}\;,
$$
where $d(\cdot,\cdot)$ denotes the distance on $N_i$ with respect to
the induced Riemannian metric.  We shall use the following estimate on
the volumes of the ``cusp'': there exists $c_2>0$ such that for every
$R>0$,
\begin{equation}\label{eq:cusp}
\overline\nu_{M_i}(N_i-(N_i)_R)=\operatorname{O}(e^{-c_2 R}).
\end{equation}
To prove this estimate, we may pass to an equivalent Riemannian metric
and to a finite index subgroup of $\Gamma\cap M_i$. This way, we reduce
the proof to the case when $M_i$ is semisimple, $\Gamma\cap M_i$ is an
nonuniform lattice in $M_i$, and the Riemannian metric on $M_i$ is
bi-invariant under a maximal compact subgroup in $M_i$. Then Equation
\eqref{eq:cusp} follows from \cite[\S 5.1]{KleMar99} (which notes that
the irreducible assumption on $\Gamma\cap M_i$ is not necessary).

Equation \eqref{eq:cusp} implies that
\begin{equation}\label{eq:eq1}
\int_{N_i-(N_i)_R}\int_{Y_{U_i}}
  f(ya^{-1}n)\;d\overline\mu_{U_i}(y)\;d\overline\nu_{M_i}(n)=
\operatorname{O}\big(e^{-c_2 R}\;\|f|_{ Y_{G_i}}\|_0\big)\;.
\end{equation}
Given $m\in M_i$ such that $(\Gamma\cap M_i)m\in (N_i)_R$, there
exists $m'\in M_i$ such that $(\Gamma\cap M_i)m=(\Gamma\cap M_i)m'$
and $d(e,m')\le R$.  Therefore, it follows from Equation
\eqref{eq:km2} that
\begin{align}
  &\int_{(N_i)_R}\int_{Y_{U_i}} f(ya^{-1}n)\;
  d\overline\mu_{U_i}(y)\;d\overline\nu_{M_i}(n) \nonumber\\
  =\;\;&\overline\nu_{M_i}((N_i)_R) (\mathcal{P}_{G_i} f)(\Gamma e)
  +\operatorname{O}\big(E_i(a)^{\kappa_1} e^{c_1\, R}\;
  \|f|_{Y_{G_i}}\|_\ellk\big)
  \nonumber\\
  =\;\;& (\mathcal{P}_{G_i} f)(\Gamma e)
  +\operatorname{O}\big((e^{-c_2 R}+ E_i(a)^{\kappa_1} e^{c_1 R})
  \;\|f|_{Y_{G_i}}\|_\ellk\big)\;.\label{eq:eq2}
\end{align}
Finally, combining \eqref{eq:eq1} and \eqref{eq:eq2}, we obtain that
\begin{align*}
  &\int_{N_i}\int_{Y_{U_i}}
  f(ya^{-1}n)\;d\overline\mu_{U_i}(y)\;d\overline\nu_{M_i}(n)
  \\=\;\;&(\mathcal{P}_{G_i} f)(\Gamma e) +\operatorname{O}\big((e^{-c_2 R}+
  E_i(a)^{\kappa_1} e^{c_1 R})\;\|f|_{Y_{G_i}}\|_\ellk\big)\;.
\end{align*}
Taking $R=\log E_i(a)^{-\frac{\kappa_1}{c_1+c_2}}$, we deduce the
claim of Lemma \ref{lem:step_H} with $\kappa_2=
\frac{\kappa_1c_2}{c_1+c_2}$.  \cqfd

\bigskip\noindent{\bf Proof of Proposition \ref{th:equidist_H}. }
For a subsemigroup $D$ which decomposes as a product $D=D_p\cdots D_q$
and $p\leq i\leq q$,
we write
\begin{align*}
D_{\le i}=D_p\cdots D_i\quad\hbox{and}\quad D_{> i}=D_{i+1}\cdots D_q.
\end{align*}
We show inductively on $i\in\{0,\dots, s\}$ that for every $a=a_1\dots
a_i\in A_{\le i}^+$ (by convention $a=e$ if $i=0$) and $g\in
G_{>i}$, we have
\begin{align}\label{eq:induction}
\int_{Y_{H_{\le i}}} f(ya^{-1}g)\, d\overline\mu_{H_{\le i}}(y)
=(\mathcal{P}_{G_{\le i}}f)(\Gamma g)
+\sum_{j=1}^i \operatorname{O}\big( E_j(a_j)^{\kappa} \| f \|_\ellk \big)
\end{align}
with $\kappa=\kappa_2$ and $\ellk$ as in Lemma \ref{lem:step_H}. 
Since $H_{\le 0}=G_{\le 0}=Z(G)$, this is obvious for $i=0$.
To get this estimate for $i=1$, we apply Lemma \ref{lem:step_H}
to the function $f_g(y)=(\mathcal{P}_{G_{\le 0}}f)(yg)$ with $g\in G_{>1}$.
Since $G_1$ commutes with $G_{\le 0}$ and $G_{>1}$, we have
$$
\| f_g|_{Y_{G_1}} \|_\ellk\le 
\| f \|_\ellk\quad\hbox{and}\quad 
(\mathcal{P}_{G_1}\mathcal{P}_{G_{\le 0}}f)(\Gamma g) =
(\mathcal{P}_{G_{\le 1}}f)(\Gamma g)\;.
$$
This proves Equation \eqref{eq:induction} with $i=1$.

Now suppose that Equation \eqref{eq:induction} is proved at rank $i$.
As in Equation \eqref{eq:Y_H}, for $f\in \C_c(Y_G)$,
$$
\int_{Y_{H_{\le i+1}}} f(y)\, d\overline\mu_{H_{\le i+1}}(y)\\
=\int_{(\Gamma\cap H_i)\backslash H_i}
\int_{Y_{H_{\le i}}} f(yh)\, d\overline\mu_{H_{\le i}}(y)
\;d\overline\nu_{H_i}(h)\;.
$$
Hence, for every $a'=a_1\dots a_{i}\in A_{\le i}^+$, $a_{i+1}\in
A_{i+1}^+$ and $g\in G_{>i+1}$, with $a=a'a_{i+1}$, by the right
invariance of $\overline{\nu}_{H_{i+1}}$ under $H_{i+1}$ and by Equation
\eqref{eq:induction}, we have
\begin{align*}
&\int_{Y_{H_{\le i+1}}} f(ya^{-1}g)\, d\overline\mu_{H_{\le i+1}}(y)\\
=&\int_{(\Gamma\cap H_{i+1})\backslash H_{i+1}}
\int_{Y_{H_{\le i}}} f(y(a')^{-1}ha_{i+1}^{-1}g)\; 
d\overline\mu_{H_{\le i}}(y) \;d\overline\nu_{H_{i+1}}(h)\\
=&\int_{(\Gamma\cap H_{i+1})\backslash H_{i+1}}
(\mathcal{P}_{G_{\le i}}f)(\Gamma ha_{i+1}^{-1}g)\; 
d\overline\nu_{H_{i+1}}(h)+
\sum_{j=1}^i \operatorname{O}\big( E_j(a_j)^\kappa \| f \|_\ellk\big)\;.
\end{align*}
Applying Lemma \ref{lem:step_H} to the functions $\overline
f_g:y\mapsto(\mathcal{P}_{G_{\le i}}f)(yg)$ on $Y_{G}$, we obtain
\begin{align*}
&\int_{(\Gamma\cap H_{i+1})\backslash H_{i+1}}
(\mathcal{P}_{G_{\le i}}f)(\Gamma ha_{i+1}^{-1}g)\; d\overline\nu_{H_{i+1}}(h)
\\
=&\; (\mathcal{P}_{G_{i+1}}\overline f_g)(\Gamma e)
+ \operatorname{O}\big( E_{i+1}(a_{i+1})^\kappa 
\;\| \overline f_g|_{Y_{G_{i+1}}} \|_\ellk\big)\\
=&\; (\mathcal{P}_{G_{\le i+1}}f)(\Gamma g)
+ \operatorname{O}\big( E_{i+1}(a_{i+1})^\kappa \;\| f \|_\ellk\big).
\end{align*}
This completes the proof of Equation \eqref{eq:induction}. Since
$$
(\mathcal{P}_{G_{\le s}}f)(\Gamma e)=\int_{Y_G} f\,d\overline\mu_G,
$$
the proposition follows.
\cqfd

\medskip\noindent{\bf Proof of Proposition \ref{prop:step2}. }
Since
\begin{align*}
\int_{Y_{P_T}} f\,d\mu_P
=\int_{A_T}\Big(\int_{Y_H} f(ya^{-1})\, d\mu_H(y)\Big)\,
\Big(\prod_{\alpha\in \Delta-I} \alpha(a)^{m_\alpha}\Big)\;d\omega_A(a)\;,
\end{align*}
it follows from Proposition \ref{th:equidist_H} and from Equation
\eqref{eq:dandemo} that
\begin{align*}
\int_{Y_{P_T}} f\,d\mu_P=\mu_P(Y_{P_T})\int_{Y_G} f\, d\overline \mu_G
+\operatorname{O}\Big(\|f\|_\ellk\int_{A_T}E^\kappa(a)
\Big(\prod_{\alpha\in \Delta-I} \alpha(a)^{m_\alpha}\Big)\;d\omega_A(a) \Big).
\end{align*}
For every $i\in\{1,\dots, s\}$ such that $\Delta_i-I\neq \emptyset$, let
$\beta\in \Delta_i-I$. For every $b_i\in A_i$, we have
$$
E_i(b_i)\le e^{-\log \beta(b_i)}.
$$
Hence, by Lemma \ref{lem:croissancePT}, we have, assuming that
$\kappa<\min_{\alpha\in\Delta-I} m_\alpha$ (which is possible),
\begin{align*}
  \int_{A_T}E_i(a_i)^\kappa \Big(\prod_{\alpha\in \Delta-I}
  \alpha(a)^{m_\alpha}\Big)\;d\omega_A(a) &\leq
  c_A\Big(\prod_{\alpha\in \Delta-I-\{\beta\}} \int_0^{t_\alpha}
  e^{m_\alpha s}\,ds\Big)
  \int_0^{t_\beta} e^{(m_\beta-\kappa)s}\,ds\\
  & =\operatorname{O}\big(\mu_P(Y_{P_T})\; e^{-\kappa t_\beta}\big)\;.
\end{align*}
Therefore, since $E^\kappa(a)=\sum_{1\leq i\leq s\;:\;\Delta_i-I\neq \emptyset}
E_i(a_i)^\kappa$, we have
$$
\frac{1}{\mu_P(Y_{P_T})}\int_{Y_{P_T}} f\,d\mu_P=\int_{Y_G} f\, d\overline\mu_G
+\operatorname{O}(e^{-\kappa \min T}\|f\|_\ellk),
$$
as required.
\cqfd

\bigskip 
\noindent{\bf Step 3. } In this last step of the proof of Theorem
\ref{theo:main}, we will diffuse the orbits of $L\cap H$ we want to
count using bump functions, and apply the equidistribution result
given by Proposition \ref{prop:step2} in Step 2 to infer our main
theorem.

Before starting this program, we rewrite the sum whose asymptotic we
want to study in a more concise way. Let $T,T'\in
[0,+\infty[^{\Delta-I}$. By transversality (see
for instance \cite[p.~22, Theo.~3.3]{Hirsch76}), the intersection
$$
Z_{[T,T']}=Y_L\cap Y_{P_{[T,T']}}
$$
is a compact Riemannian submanifold of $Y_G$, invariant under the
right action of $L\cap H$, and for every $x\in Z_{[T,T']}$, we have
$T_xZ_{[T,T']} =(T_xY_L)\cap (T_xY_{P_{[T,T']}})$.  Since $\lll\cap
\ppp=\lll\cap\hhh$ by Equation \eqref{eq:lcappeglcaph}, the Lie group
$L\cap H$ has open orbits in $Z_{[T,T']}$.  Hence the compact subset
$Z_{[T,T']}$ is a finite union of orbits of $L\cap H$ (see the picture
below when $A$ is $1$-dimensional).

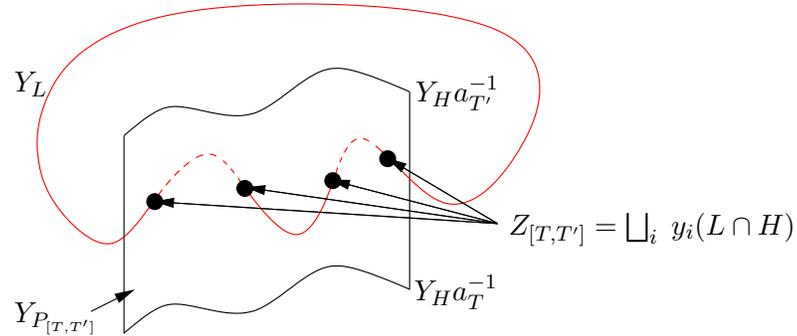
\begin{center}
\begin{picture}(0,0)%
\includegraphics{fig_intersection.pstex}%
\end{picture}%
\setlength{\unitlength}{3631sp}%
\begingroup\makeatletter\ifx\SetFigFont\undefined%
\gdef\SetFigFont#1#2#3#4#5{%
  \reset@font\fontsize{#1}{#2pt}%
  \fontfamily{#3}\fontseries{#4}\fontshape{#5}%
  \selectfont}%
\fi\endgroup%
\begin{picture}(3615,2275)(586,-1873)
\put(3331,-1636){\makebox(0,0)[lb]{\smash{{\SetFigFont{11}{13.2}{\rmdefault}{\mddefault}{\updefault}{\color[rgb]{0,0,0}$Y_Ha_{T}^{-1}$}%
}}}}
\put(3339,-278){\makebox(0,0)[lb]{\smash{{\SetFigFont{11}{13.2}{\rmdefault}{\mddefault}{\updefault}{\color[rgb]{0,0,0}$Y_Ha_{T'}^{-1}$}%
}}}}
\put(601,-1786){\makebox(0,0)[lb]{\smash{{\SetFigFont{11}{13.2}{\rmdefault}{\mddefault}{\updefault}{\color[rgb]{0,0,0}$Y_ {P_{[T,T']}}$}%
}}}}
\put(601,-211){\makebox(0,0)[lb]{\smash{{\SetFigFont{11}{13.2}{\rmdefault}{\mddefault}{\updefault}{\color[rgb]{0,0,0}$Y_L$}%
}}}}
\put(3976,-1186){\makebox(0,0)[lb]{\smash{{\SetFigFont{11}{13.2}{\rmdefault}{\mddefault}{\updefault}{\color[rgb]{0,0,0}$Z_{[T,T']}=\bigsqcup_i\;y_i(L\cap H)$}%
}}}}
\end{picture}%

\end{center}

We will denote by $\mu_{Z_{[T,T']}}$ the Riemannian measure on
$Z_{[T,T']}$. Using Riemannian volumes, we hence have
\begin{align*}
\mu_{Z_{[T,T']}}(Z_{[T,T']})&=\sum_{[y]\in (Y_L\cap Y_{P_{[T,T']}})/(L\cap H)} \;
\vol\big(y(L\cap H)\big)\\ &=
\sum_{a\in A_{[T,T']}}\;\;\sum_{[y]\in (Y_L\cap Y_H a^{-1})/(L\cap H)} \;
\vol\big(y(L\cap H)\big)\;.
\end{align*}
By Lemma \ref{lem:grouptheoretic} in Step 1, the quantity
$\mu_{Z_{[0,T]}}(Z_{[0,T]})$, when divided by $\vol(Y_L)$, is the
sum whose asymptotic we want to study.

\medskip We first start by studying the supports of the bump functions
we will define: they will be appropriate neighbourhoods of $Y_L$ and
$Z_{[T,T']}$. Fix $\epsilon>0$, which will be appropriately choosen
small enough later on. Consider the open ball $B(0,\epsilon)$ of
center $0$ and radius $\epsilon$ in the orthogonal complement
$\qqq\oplus\aaa$ of $\lll\cap\ppp$ in $\ppp$, and let
$\OOO_\epsilon=\exp B(0,\epsilon)$, which is contained in $P$.

Since $L$ is compact, if $\epsilon$ is small enough, the right action
of $G$ on $Y_G$ induces a map $Y_L\times\OOO_\epsilon\ra Y_G$, with
$(y,g)\mapsto yg$, which is a smooth diffeomorphism onto an open
neighbourhood $Y_L\OOO_\epsilon$ of the submanifold $Y_L$ in
$Y_G$. Similarly, if $\epsilon$ is small enough, then for every
$T,T'\in [0,+\infty[^{\Delta-I}$, the map $Z_{[T,T']} \times
\OOO_\epsilon \ra Y_P$ defined by $(y,g)\mapsto yg$ is a smooth
diffeomorphism onto an open neighbourhood $Z_{[T,T']}\OOO_\epsilon$ of
the submanifold $Z_{[T,T']}$ in $Y_P$. If $\eta\in\RR$ and
$T''=(t''_\alpha)_{\alpha\in\Delta-I}\in [0,+\infty[^{\Delta-I}$, we
denote $T''+\eta= (t''_\alpha+\eta)_{\alpha\in\Delta-I}$.

\blemm \label{lem:doubleinclusion} 
There exists $c>0$ such that if $\epsilon>0$ is small enough, for
every $T,T'\in [0,+\infty[^{\Delta-I}$, then
$$
Z_{[T+c\epsilon,T'-c\epsilon]}\OOO_\epsilon\subset
Y_L\OOO_\epsilon\cap Y_{P_{[T,T']}}\subset 
Z_{[T-c\epsilon,T'+c\epsilon]}\OOO_\epsilon\;.
$$
\elemm

\dem
We first claim that there exists $c>0$ such that 
$$
P_{[T+c\epsilon,T'-c\epsilon]}\subset P_{[T,T']}\OOO_\epsilon
\subset P_{[T-c\epsilon,T'+c\epsilon]}\;.
$$ 
Since the product map $(h,a)\mapsto ha$ is a diffeomorphism from
$H\times A$ to $P$, since $\OOO_\epsilon$ is contained in $P$, and
since the distances are Riemannian ones, there exists $c_1>0$ such
that if $\epsilon>0$ is small enough, then for every $g\in
\OOO_\epsilon$, there exist $h\in H$ and $a\in A$ with $g=ha$ and
$d(a,e)\leq c_1\epsilon$. Since the Riemannian distance on $A$ is
equivalent to the image by $\exp$ of the distance on $\aaa$ defined by
the norm $\|x\|=\max_{\alpha\in\Delta-I}|\log (\alpha (\exp x))|$,
there exists $c_2>0$ such that $|\log \alpha (a)|\leq c_2d(a,e)$ for
every $a\in A$.

Let $g\in \OOO_\epsilon$, $h\in H$ and $a\in A$ be such that $g=ha$
and $d(a,e)\leq c_1\epsilon$. Since $A$ normalises $H$, we have
$HA^{-1}_{[T,T']}g=HA^{-1}_{[T,T']}ha=HA^{-1}_{[T,T']}a $. Hence
$HA^{-1}_{[T,T']}g$ is contained in $HA^{-1}_{[T-c_1c_2\epsilon,
  T'+c_1c_2\epsilon]}$ and contains $HA^{-1}_{[T+c_1c_2\epsilon,
  T'-c_1c_2\epsilon]}$. This proves the first claim.  

Now, let $y\in Y_L$, $g\in \OOO_\epsilon$ and $p\in P_{[T,T']}$ be
such that $yg=\pi(p)$. Then $y=\pi(pg^{-1})$. Since $\OOO_\epsilon$ is
invariant by taking inverses, $pg^{-1}$ belongs to
$P_{[T,T']}\OOO_\epsilon$, hence by the first claim, $yg\in
Z_{[T-c\epsilon,T'+c\epsilon]} \OOO_\epsilon$. The left inclusion is
proven similarly.  
\cqfd

\medskip We now study the properties of the Riemannian measures on
the neighbourhoods $Y_L\OOO_\epsilon$ and $Z_{[T,T']}\OOO_\epsilon$.

\blemm \label{lem:desintegr} For every $\epsilon>0$ small enough,
there exist smooth measures $\nu$ and $\wt\nu$ on $\OOO_\epsilon$ such
that the product maps $Y_L\times\OOO_\epsilon\ra Y_G$ and $Z_{[T,T']}
\times \OOO_\epsilon \ra Y_P$ send the product measures $\mu_L\otimes
\nu$ and $\mu_{Z_{[T,T']}}\otimes \wt \nu$ to the restricted measures
${\mu_G}_{\mid Y_L\OOO_\epsilon}$ and ${\mu_P}_{\mid Z_{[T,T']}
  \OOO_\epsilon}$, respectively.  Furthermore, $\frac{d\wt \nu}
{d\nu}(e)=1$.  \elemm

\dem Since the measure ${\mu_G}_{\mid Y_L\OOO_\epsilon}$ (respectively
${\mu_P}_{\mid Z_{[T,T']} \OOO_\epsilon}$) is Riemannian, it
disintegrates with respect to the trivialisable fibration
$Y_L\OOO_\epsilon\ra Y_L$ (respectively $Z_{[T,T']}\OOO_\epsilon\ra
Z_{[T,T']}$ with measure on the basis $\mu_L$ (respectively
$\mu_{Z_{[T,T']}}$), and conditional measures $\nu_y$ (respectively
$\wt \nu_y$ on the fibers $y\OOO_\epsilon$ for all $y\in Y_L$
(respectively $y\in Z_{[T,T']}$). By left invariance of the measures
$\omega_L$ and $\omega_{L\cap H}$, there exist smooth measures $\nu$
(respectively $\wt \nu$) on $\OOO_\epsilon$ such that the maps
$\OOO_\epsilon\ra y\OOO_\epsilon$ defined by $g\mapsto yg$ send $\nu$
(respectively $\wt \nu$) to $\nu_y$ (respectively $\wt \nu_y$) for all
$y\in Y_L$ (respectively $y\in Z_{[T,T']}$). This proves the first
claim.

Since $\qqq+\aaa$ is orthogonal to $\lll$ (respectively
$\lll\cap\hhh$) by Equation \eqref{eq:perpendicul}, the manifold
$y\OOO_\epsilon$ is orthogonal to $Y_L$ (respectively $Z_{[T,T']}$) at
every $y\in Y_L$ (respectively $y\in Z_{[T,T']}$). Hence any
orthonormal frame $F$ of $T_y (y\OOO_\epsilon)$ at a given $y\in
Z_{[T,T']}$ may be completed to an orthonormal frame whose last
vectors form a basis of $T_y Y_L$, whose first vectors form a basis of
$T_y Z_{[T,T']}$. By desintegration, the orthogonal frame $F$ has
the same infinitesimal volume for $\nu$ and $\wt \nu$. The last
assertion follows.  \cqfd

\medskip Let us now define our bump functions. By the standard
construction of bump functions on manifolds, for every $\ellk\in\NN$,
there exists $\kappa'>0$ such that for every $\epsilon>0$ small
enough, there exists a ${\rm C}^\ellk$ map $\psi_\epsilon$ from
$\OOO_\epsilon$ to $[0,+\infty[\,$, with compact support, such that
$\int \psi_\epsilon\;d\nu=1$ and $\|\psi_\epsilon\|_\ellk=
\operatorname{O}(\epsilon^{-\kappa'})$.  Since $\frac{d\wt \nu}
{d\nu}=1+\operatorname{O}(\epsilon)$ on $\OOO_\epsilon$ by Lemma
\ref{lem:desintegr}, we have
$$
\int_{\OOO_\epsilon} \psi_\epsilon\;d\wt \nu=1+\operatorname{O}(\epsilon)\;.
$$

For every $\epsilon>0$ small enough, define $f_\epsilon:Y_G \ra
[0,+\infty[$ by $f_\epsilon(y)=0$ if $y\notin Y_L\OOO_\epsilon$ and
$f_\epsilon(yg)=\psi_\epsilon(g)$ for every $y\in Y_L$ and $g\in
\OOO_\epsilon$.  Note that $f_\epsilon$ is ${\rm C}^\ellk$ with compact
support, since $Y_L$ is compact. We have 
$$
\int_{Y_G} f_\epsilon\;d\overline\mu=\frac{\int_{Y_L\OOO_\epsilon}
  f_\epsilon\;d\mu_G}{\vol(Y_G)}= \frac{\int_{g\in \OOO_\epsilon}
\int_{y\in Y_L}\psi_\epsilon(g)\;d\mu_L(y)d\nu(g)}{\vol(Y_G)}
=\frac{\vol(Y_L)}{\vol(Y_G)}\;,
$$
and $\|f_\epsilon\|_\ellk= \operatorname{O}(\epsilon^{-\kappa'})$.

Since the support of $f_\epsilon$ is contained in $Y_L\OOO_\epsilon$,
by Lemma \ref{lem:desintegr}, and by the right inclusion in Lemma
\ref{lem:doubleinclusion}, we have, for every $T\in
[0,+\infty[^{\Delta-I}$,
\begin{align}
  \int_{Y_{P_T}} f_\epsilon \;d\mu_P&\leq
  \int_{Z_{[-c\epsilon,T+c\epsilon]}\OOO_\epsilon} f_\epsilon
  \;d\mu_P\nonumber\\ &= \int_{g\in \OOO_\epsilon}\int_{y\in
    Z_{[-c\epsilon,T+c\epsilon]}} \psi_\epsilon(g)
  \;d\mu_{Z_{[-c\epsilon,T+c\epsilon]}}(y)d\wt\nu(g)\nonumber\\ &\label{eq:majo}
  =\vol(Z_{[-c\epsilon,T+c\epsilon]})\big(1+\operatorname{O}(\epsilon)\big)\;.
\end{align}
Similarly, since $f_\epsilon\geq 0$ and by the left inclusion in Lemma
\ref{lem:doubleinclusion}, we have, for every $T\in
[0,+\infty[^{\Delta-I}$,
\begin{align}\label{eq:mino}
  \int_{Y_{P_T}} f_\epsilon \;d\mu_P\geq
  \int_{Z_{[c\epsilon,T-c\epsilon]}\OOO_\epsilon} f_\epsilon
  \;d\mu_P
  =\vol(Z_{[c\epsilon,T-c\epsilon]})\big(1+\operatorname{O}(\epsilon)\big)\;.
\end{align}

Finally, we apply Step 2 to our bump functions.  By Proposition
\ref{prop:step2}, we have the equality $ \frac{1}{\mu_P(Y_{P_T})}
\int_{Y_{P_T}} f_\epsilon \;d\mu_P\;= \int_{Y_G} f_\epsilon
\;d\overline{\mu}_G+ \operatorname{O}\big( \;e^{-\kappa\min
  T}\,\|f_\epsilon\|_\ellk\big) $. Hence, by the properties of
$f_\epsilon$,
\begin{align}\label{eq:applistep2}
\int_{Y_{P_T}}f_\epsilon\;d\mu_P=
\frac{\vol(Y_L)\mu_P(Y_{P_T})}{\vol(Y_G)}\big(1
+\operatorname{O}(\epsilon^{-\kappa'}e^{-\kappa \min T})\big)\;.
\end{align}

Let $\delta=\frac{\kappa}{\kappa'+1}>0$ and $\epsilon= e^{-\delta\min
  T}$ (which tends to $0$ as $\min T$ tends to $+\infty$). Then
$\epsilon^{-\kappa'}e^{-\kappa \min T} = e^{(\kappa'\delta-\kappa)\min
  T} =e^{-\delta\min T}$.  By the equations \eqref{eq:mino} and
\eqref{eq:applistep2}, and by Lemma \ref{lem:croissancePT}, we have,
as $\min T$ tends to $+\infty$,
\begin{align*}
  \vol(Z_{[c\epsilon,T-c\epsilon]})&\leq \Big(\int_{Y_{P_T}}
  f_\epsilon \;d\mu_P\Big)\big(1+\operatorname{O}(e^{-\delta\min
    T})\big)\\& =\frac{\vol(Y_L)\mu_P(Y_{P_T})}{\vol(Y_G)}
  \big(1+\operatorname{O}(e^{-\delta\min T})\big)\\
  &=\frac{\Vol(\Lambda^{\vee}\bs
    A)\vol(Y_L)\vol(Y_H)}{\vol(Y_G)}\Big(\prod_{\alpha\in
    \Delta-I}\frac{e^{m_\alpha
      t_\alpha}}{m_\alpha}\Big)\big(1+\operatorname{O}(e^{-\delta\min
    T})\big) \;.
\end{align*}

Since $e^x=1+\operatorname{O}(x)$ as $x$ tends to $0$, we have $e^{c\,
  e^{-\delta \min T} \sum_{\alpha\in \Delta-I}m_\alpha }= 1+
\operatorname{O}(e^{-\delta\min T})$ as $\min T$ tends to $+\infty$.
Since $Z_{[0,c\epsilon]}$ is bounded, we hence have, as $\min T$ tends
to $+\infty$,
$$
\vol(Z_{[0,T]})\leq
\frac{\Vol(\Lambda^{\vee}\bs
  A)\vol(Y_L)\vol(Y_H)}{\vol(Y_G)}\Big(\prod_{\alpha\in
  \Delta-I}\frac{e^{m_\alpha
    t_\alpha}}{m_\alpha}\Big)\big(1+\operatorname{O}(e^{-\delta\min
  T})\big)\;.
$$
The converse inequality is proven similarly, using Equation
\eqref{eq:majo} instead of Equation \eqref{eq:mino}

\medskip 
Since $\sum_{a\in A_{T}}\;\; \sum_{[x]\in (L\cap \Ga) \bs
  (Lv_a\cap \Ga v_0)} \; {w'}_{{\bf L},\rho_{\mid {\bf L}}}(x) =
\frac{\vol(Z_{[0,T]})}{\vol(Y_L)}$ as said in the beginning of Step 3,
this ends the proof of Theorem \ref{theo:main}.  
\cqfd

\medskip\noindent\brema\label{rem:stansiegelweight}
{\rm Let ${\bf G,P,A,M,U,L,V}, \rho, v_0$ be as in
the statement of Theorem \ref{theo:main}, and assume furthermore that
${\bf G}$ is simply connected. Then we have the following counting
results using the standard Siegel weights.

There exists $\delta>0$ such that, as $T=
(t_\alpha)_{\alpha\in\Delta-I} \in[0,+\infty[^{\Delta-I}$ and
$\min_{\alpha\in\Delta- I} t_\alpha$ tends to $+\infty$,
$$
\sum_{a\in A_T}\;\;\;
\sum_{[x]\in {\bf L}(\ZZ)\bs(\rho({\bf
    L}(\RR) a)v_0\cap \rho(\Ga)v_0)} \; 
w_{{\bf L},\rho_{\mid {\bf L}}}(x)= 
$$
$$
\frac{\vol\big({\bf MU}(\ZZ)\bs {\bf MU}(\ZZ)\big)
\vol(\Lambda^{\!\vee}\,\bs {\bf A}(\RR)_0)}
{\vol({\bf G}(\ZZ)\bs {\bf G}(\RR))}\;
\Big(\prod_{\alpha\in \Delta- I}
\frac{e^{m_\alpha t_\alpha}}{m_\alpha}\Big)\big(1+
\operatorname{O}(e^{-\delta\min_{\alpha\in\Delta- I}t_\alpha})\big)\;.
$$}
\erema

The proof is the same as the one of Theorem \ref{theo:main}, with the
following modifications. Since $\bf G$ is simply connected, ${\bf
  G}(\RR)$ is connected (see for instance \cite[\S 7.2]{PlaRap94}).
Hence with the previous notation, we have $G={\bf G}(\RR)$ and
$\Ga=\Ga(\ZZ)$ (and the connectedness of $G$ was useful). Now take
$L={\bf L}(\RR)$ instead of $L={\bf L}(\RR)_0$ (which is still
contained in $G$, but would not have been if $G$ was only taken to be
${\bf G}(\RR)_0$ while ${\bf G}(\RR)$ is not connected). Though $L$ and
$Y_L$ may be no longer connected, the proof stays valid.

\bigskip 
To end this section, we give two slightly different versions
of Theorem \ref{theo:main} when $\bf P$ is maximal.

\btheo \label{theo:mainmaximal} Let $\bf G$ be a connected reductive
linear algebraic group defined over $\QQ$, without nontrivial
$\QQ$-characters.  Let $\bf P$ be a maximal (proper) parabolic
subgroup of $\bf G$ defined over $\QQ$, and let ${\bf P} ={\bf AMU}$
be a relative Langlands decomposition of $\bf P$, such that ${\bf A}
(\RR)_0$ is a one-parameter subgroup $(a_s)_{s\in\RR}$, with
$\lambda=\log \det\, (\Ad a_1)_{\mid \UUU}>0$, where $\UUU$ is the Lie
algebra of ${\bf U} (\RR)$.  Let $\rho:{\bf G} \ra \GL({\bf V})$ be a
rational representation of $\bf G$ defined over $\QQ$ such that there
exists $v_0\in{\bf V}(\QQ)$ whose stabiliser in $\bf G$ is $\bf
MU$. Let $\bf L$ be a reductive algebraic subgroup of $\bf G$ defined
and anisotropic over $\QQ$. Assume that $\bf LP$ is Zariski-open in
$\bf G$ and that for every $s\in\RR$, the orbit ${\bf X}_s=\rho({\bf
  L} a_s)v_0$ is Zariski-closed in $\bf V$.

\medskip (1) Endow ${\bf G} (\RR)$ with a left-invariant Riemannian
metric, for which the Lie algebras of ${\bf MU}(\RR)$ and ${\bf
  A}(\RR)$ are orthogonal, and the orthogonal of the Lie algebra of
${\bf P}(\RR)$ is contained in the Lie algebra of ${\bf L}(\RR)$.  Let
$G={\bf G}(\RR)_0$ and $\Ga={\bf G}(\ZZ)\cap G$.  There exists
$\delta>0$ such that, as $t\geq 0$ tends to $+\infty$,
\begin{align*}
\sum_{0\leq s\leq t}\;\;\;&\sum_{[x]\in 
({\bf L}(\RR)_0\cap\Ga)\bs
(\rho({\bf L}(\RR)_0 a_s)v_0\cap \rho(\Ga)v_0)} \;
{w'}_{{\bf L},\rho_{\mid {\bf L}}}(x)\\ &\;\;\;\;\;\;\;\;=\;\; 
\frac{\vol\big(({\bf MU}\cap\Ga)\bs ({\bf MU}\cap G)\big)
\vol(a_1^\ZZ\bs {\bf A}(\RR)_0)}{\lambda\vol(\Ga\bs G)}\;e^{\lambda t}+
\operatorname{O}(e^{(\lambda-\delta)t})\;.
\end{align*}

(2) Let $\Lambda$ be a $\ZZ$-lattice in ${\bf V}(\QQ)$ invariant under
${\bf G}(\ZZ)$, and let $\Lambda^{\rm prim}$ be the subset of
indivisible elements of $\Lambda$. Assume $\rho$ to be irreducible
over $\CC$. Then there exist $c,\delta>0$ such that, as $t\geq 0$
tends to $+\infty$,
$$
\sum_{0\leq s\leq t}\;\;\;\sum_{[x]\in 
({\bf L}(\ZZ)\cap{\bf L}(\RR)_0)\bs({\bf X}_s\cap \Lambda^{\rm prim})} \;
w'_{{\bf L},\rho_{\mid {\bf L}}}(x)= c\;e^{\lambda t}+
\operatorname{O}(e^{(\lambda-\delta)t})\;.
$$
\etheo

\dem (1) In this case, $\Delta-I$ consists of one simple root
$\alpha_0$. Changing the parametrisation of the one-parameter subgroup
$(a_s)_{s\in\RR}$ appearing in Theorem \ref{theo:mainmaximal} by
multiplying $s$ by a positive constant does not change the asymptotic
formula in the statement of Theorem \ref{theo:mainmaximal} (1). Hence
we may assume that $a_1=(\alpha_0)^\vee$, hence that the group
$a_1^\ZZ$ generated by $a_1$ is equal to the lattice $\Lambda^\vee$.
The constant $\lambda$ defined in Theorem \ref{theo:mainmaximal} is
then equal to $m_{\alpha_0}$. The first part of Theorem
\ref{theo:mainmaximal} hence follows from Theorem \ref{theo:main}.

\medskip (2)
We start by proving two lemmas.

\blemm  \label{lem:irreducpluopoi}
If $\rho$ is irreducible, then the stabiliser of $\CC v_0$ in
$\bf G$ is $\bf P$ and there exists $\chi\in\RR$ such that
$a_sv_0=e^{\chi s}v_0$ for every $s\in\RR$.  
\elemm

\dem Let ${\bf T}$ be a maximal torus of $\bf G$ containing $\bf S$,
and let $\Delta_{\bf T}$ be a set of primitive roots of $\bf G$
relative to $\bf T$, whose set of nonzero restrictions to $\bf S$ is
$\Delta$ (see for instance \cite[\S 21.8]{Borel91}. Then the unipotent
subgroup ${\bf U}^+_{\bf T}$, whose Lie algebra is the sum of the
positive root spaces of $\bf G$ relative to $\bf T$, is contained in
$\bf MU$. By the properties of the highest weights, if $\rho$ is
irreducible, the space $\{v\in{\bf V}\;:\;{\bf U}^+_{\bf T} v=v\}$ is
one-dimensional, hence equal to $\CC v_0$. Since $\bf A$ normalises
$\bf MU$, hence ${\bf U}^+_{\bf T}$, it preserves $\CC v_0$, and the
result follows, by the connectedness of $A$. \cqfd

\blemm \label{lem:decompXprimorb}
There exist $v_1,\dots, v_k$ in $\Lambda^{prim}$ such that
$\Lambda^{\rm prim}\cap{\bf G}v_0=\bigsqcup_{i=1}^k\Ga v_i$.  
\elemm

\dem By \cite[Prop.~20.5]{Borel91}, the natural map ${\bf
  G}(\mathbb{Q})\to ({\bf G}/{\bf P})(\mathbb{Q})$ is onto. Since
${\bf G}v_0\simeq {\bf G}/{\bf MU}$, this implies that every $x\in
({\bf G}v_0)(\mathbb{Q})$ may be written as $x=gpv_0$ for some $g\in {\bf
  G}(\mathbb{Q})$ and $p\in {\bf P}$. Hence by Lemma
\ref{lem:irreducpluopoi},
$$
({\bf G}v_0)(\mathbb{Q})\subset \mathbb{C}^\times {\bf G}(\mathbb{Q})v_0.
$$
By \cite[Prop.~15.6]{Borel69}, there exists a finite subset $F$ of
${\bf G}(\mathbb{Q})$ such that ${\bf G}(\mathbb{Q})=\Gamma F {\bf
  P}(\mathbb{Q})$. Hence,
$$
({\bf G}v_0)(\mathbb{Q})\subset \mathbb{C}^\times \Gamma F v_0.
$$
In particular, we conclude that there exist $v_1,\ldots v_k$ in
$\Lambda^{prim}$ such that
$$
\Lambda^{prim}\cap {\bf G}v_0\subset \bigsqcup_{i=1}^k
\mathbb{C}^\times \Gamma v_i.
$$
Since for every $v\in \Lambda^{prim}$,
$$
\mathbb{C}^\times v\cap \Lambda^{prim}=\{\pm v\},
$$
this implies the lemma.\cqfd

\medskip Now, since the identity component $L$ of ${\bf L}(\RR)$ has
finite index in ${\bf L}(\RR)$, there exist $\ell_1\dots, \ell_{k'}$ in
${\bf L}(\RR)$ such that ${\bf L}(\RR)= \bigsqcup_{j=1}^{k'} L \, \ell_j$.
Hence, since $v_0$ belongs to ${\bf V}(\RR)$ and ${\bf X}_s \subset
{\bf G}v_0$, by Lemma \ref{lem:irreducpluopoi} and Lemma
\ref{lem:decompXprimorb}, we have
\begin{equation}\label{eq:trucenij}
{\bf X}_s\cap\Lambda^{\rm prim}
=({\bf L}(\RR)e^{\chi s}v_0)\cap(\Lambda^{\rm prim}\cap{\bf G}v_0)
=\bigsqcup_{1\leq i\leq k\;,\;1\leq j\leq k'}
\;e^{\chi s}L\,\ell_j v_0\cap\Ga v_i\;.
\end{equation}
If $L\,\ell_j v_0\cap\Ga v_i$ is nonempty, fix $v_{i,j}\in L\ell_j v_0
\cap \Ga v_i$. In particular, there exist $\ga\in\Ga$ and $\ell\in L$
such that $v_{i,j}=\ell\, \ell_jv_0=\ga v_i$. Since $v_i\in {\bf V}
(\QQ)$, we have $v_{i,j}\in {\bf V}(\QQ)$. Hence the stabiliser ${\bf
  P}_{i,j}$ of $v_{i,j}$ in $\bf G$ is an algebraic subgroup defined
over $\QQ$. Since $v_{i,j}$ is in the $\bf G$-orbit of $v_0$, the
stabilisers of $v_0$ and of $v_{i,j}$ are conjugate, hence ${\bf
  P}_{i,j}$ is a parabolic subgroup of $\bf G$. Since two parabolic
subgroups of $\bf G$, which are defined over $\QQ$ and conjugate in
$\bf G$, are conjugated by an element of ${\bf G}(\QQ)$ (see for
instance \cite[Theo.~20.9 (iii)]{Borel66}), there exists
$\alpha_{i,j}\in{\bf G}(\QQ)$ such that ${\bf P}_{i,j}=
\alpha_{i,j}{\bf P} \alpha_{i,j}^{-1}$.  Furthermore, using Lemma
\ref{lem:irreducpluopoi}, we have $\CC v_{i,j}=\CC\alpha_{i,j}v_0$. A
relative Langlands decomposition of ${\bf P}_{i,j}$ is ${\bf P} _{i,j}
= {\bf A}_{i,j}{\bf M}_{i,j}{\bf U}_{i,j}$ where
$$
{\bf A} _{i,j}
=\alpha_{i,j}{\bf P} \alpha_{i,j}^{-1},\;\;\;{\bf M}_{i,j}=
\alpha_{i,j}{\bf M} \alpha_{i,j}^{-1},\;\;\;{\bf U}_{i,j}=
\alpha_{i,j}{\bf U} \alpha_{i,j}^{-1}\;.
$$
We have ${\bf A}_{i,j}(\RR)_0 = \big(a_s^{i,j}=\alpha_{i,j}a_s
\alpha_{i,j}^{-1}\big)_{s\in\RR}$ and the Lie algebra of ${\bf
  U}_{i,j}(\RR)$ is $\UUU_{i,j}=\Ad \alpha_{i,j}(\UUU)$. Hence
$a_s^{i,j}v_{i,j}=e^{\chi s}v_{i,j}$ for every $s\in\RR$ and
$$
\log \det (\Ad a_1^{i,j})_{\mid \UUU_{i,j}}=\lambda\;,
$$
for every $i,j$ with $L\,\ell_j v_0\cap\Ga v_i\neq\emptyset$.

By Assertion (1) of Theorem \ref{theo:mainmaximal} applied to the
(maximal) parabolic subgroup ${\bf P}_{i,j}$ defined over $\QQ$, there
exist $c_{i,j},\delta_{i,j}>0$ (with $c_{i,j}$ explicit) such that, as
$t\geq 0$ tends to $+\infty$,
$$
\sum_{0\leq s\leq t}\;\;\;\sum_{[x]\in 
(L\cap\Ga)\bs(La_s^{i,j}v_{i,j}\cap \Ga v_{i,j})} \;
{w'}_{{\bf L},\rho_{\mid {\bf L}}}(x)= c_{i,j}\;e^{\lambda t}+
\operatorname{O}(e^{(\lambda-\delta_{i,j})t})\;.
$$
Hence, using the equations \eqref{eq:interintersimpl} and
\eqref{eq:trucenij}, with $\delta=\min_{i,j}\delta_{i,j}$ and
$c=\sum_{i,j}c_{i,j}$, we have, as $t\geq 0$ tends to $+\infty$,
\begin{align*}&\sum_{0\leq s\leq t}\;\;\;\sum_{[x]\in 
({\bf L}(\ZZ)\cap{\bf L}(\RR)_0)\bs({\bf X}_s\cap \Lambda^{\rm prim})} \;
{w'}_{{\bf L},\rho_{\mid {\bf L}}}(x)\\ =\;\;\; &
\sum_{\tiny\begin{array}{c}
1\leq i\leq k\\1\leq j\leq k'\\L\ell_j v_0\cap\Ga v_i\neq \emptyset
\end{array}}\sum_{0\leq s\leq t}\;\;\;\sum_{[x]\in 
(L\cap\Ga)\bs(La_s^{i,j}v_{i,j}\cap \Ga v_{i,j})} \;
{w'}_{{\bf L},\rho_{\mid {\bf L}}}(x)\\ =\;\;\; & c\;e^{\lambda t}+
\operatorname{O}(e^{(\lambda-\delta)t})\;.
\end{align*}
This ends the proof of Assertion (2) of Theorem
\ref{theo:mainmaximal}.  
\cqfd

\medskip\rem Using Remark \ref{rem:stansiegelweight} instead of
Theorem \ref{theo:main} in the above proof gives Theorem
\ref{theo:intromainsc} and Theorem \ref{theo:intromaincourtsc} in the
introduction.

\section{Examples}
\label{sect:applications}

In this section, we give several examples to illustrate
our main results.

\subsection{Counting inequivalent representations of integers by 
 norm forms}
\label{sect:normform}

Let $n\geq 2$. A {\it decomposable form} $F(x_1,\dots, x_n)$ is a
polynomial in $n$ variables with coefficients in $\QQ$ which is the
product of $d$ linear forms with coefficients in $\overline{\QQ}$. In
particular, a {\it norm form} is a decomposable form
$N_{K/\QQ}(\alpha_1x_1 +\dots + \alpha_n x_n)$ where $\alpha_1, \dots,
\alpha_n$ are fixed elements in a number field $K$ of degree $d$ and
$x_1,\dots, x_n$ are rational variables.  Existence of integral
solutions to the equation $F(x)=m$ is a difficult problem (see, for
instance, \cite{ColXu09}).  Here we demonstrate how our main result
applies to integral solutions of the inequality $|F(x)|\leq m$, which
can be also studied using elementary methods as in
\cite[Ch.~VI]{Lang94}. See also, following an approach of Linnik and
Sarnak, the papers \cite{EskOh06,GanOh03,Oh04}, using dilations of
relatively compact subsets.

\bcoro \label{theo:intronormform} Let $n\geq 2$, let
$F\in\QQ[x_1,\dots, x_n]$ be a rational polynomial in $n$ variables,
which is irreducible over $\QQ$, splits as a product of $n$ linearly
independant over $\CC$ linear forms with coefficients in
$\overline{\QQ}$, and satisfies $F^{-1}(]0,+\infty[)\neq \emptyset$.
Let $\Ga_F= \{g\in\SL_n(\ZZ) \;:\;F\circ g=F\}$, and for every
$k\in\QQ$, let $\Sigma_k$ be the set of $x\in\ZZ^n$ such that
$F(x)=k$.  Then there exist $c=c(F)>0$ and $\delta= \delta(n) \in\;
]0,1[$ such that, as $r\ra+\infty$,
$$
\sum_{k\in[1,r]} \; \card\big(\Ga_F\bs \Sigma_k\big)=
c\;r+\operatorname{O}\big(r^\delta\big)\;.
$$
\ecoro

It is easy to see that the above sum is indeed finite, and that the
irreducibility assumption is necessary. With $\bf L$ the stabiliser of
$F$ in $\SL_n(\CC)$, ${\bf V}=\CC^n$, $\Lambda = \ZZ^n$ and $\pi$ the
inclusion of ${\bf L}$ in $\GL({\bf V})$, this result fits into the
program described in the beginning of the introduction, since
$\Ga_F={\bf L}(\ZZ)$, the algebraic torus $\bf L$ is anisotropic over
$\QQ$ (see Lemma \ref{lem:anisotropenormform}) and acts simply
transitively on the affine subvariety of $\bf V$ with equation
$F(x)=k$ if $k\neq 0$, noting that the Siegel weights $w_{{\bf
    L},\pi}(u)=1/\vol\big( {\bf L} (\ZZ)\bs {\bf L}(\RR)\big)$ are
then constant.

When $K$ is a number field of degree $n$ with ring of integers
$\OOO_K$, taking an integral basis $(\alpha_1,\dots,\alpha_n)$ of $K$,
and $F(x_1,\dots, x_n)$ the particular norm form
$N_{K/\QQ}(\alpha_1x_1 +\dots+ \alpha_n x_n)$, we recover the
well-known counting result of the number of nonzero integral ideals
of $\OOO_K$ with trivial ideal class and norm at most $s$ (see for
instance \cite[Theorem 3, page 132]{Lang94}), giving
\begin{equation}\label{eq:normformstand}
  \{\aaa{\rm ~ideal~in~}\OOO_K\;:\;N_{K/\QQ}(\aaa)\leq s\}=
  \frac{2^{r_1}(2\pi)^{r_2}R_Kh_K}{\omega_K\sqrt{|D_K|}}\;s+
\operatorname{O}(s^{1-\epsilon})\;,
\end{equation}
where $r_1$ and $r_2$ are the numbers of real and complex conjugate
embeddings of $K$, $R_K$ is the regulator of $K$, $h_K$ is the ideal
class number of $K$, $\omega_K$ is the number of roots of unity of
$\OOO_K$, $D_K$ is the discriminant of $K$ and $\epsilon=1/n$.

\bigskip
\noindent{\bf Proof of Corollary \ref{theo:intronormform}. } 
In order to apply Theorem \ref{theo:intromaincourtsc}, let us first
define the objects appearing in its statement.

Let ${\bf G}=\SL_n(\CC)$ which is a ($\QQ$-split) quasi-simple simply
connected linear algebraic group without nontrivial
$\QQ$-characters. Let ${\bf V}=\CC^n$, $\Lambda=\ZZ^n$ (which is a
$\ZZ$-lattice in ${\bf V}(\QQ)$ invariant under ${\bf G}(\ZZ)$),
$(e_1,\dots ,e_n)$ the canonical basis of $\bf V$ and $\rho:{\bf G}\ra
\GL({\bf V})$ the monomorphism mapping a matrix $x$ to the linear
automorphism of $\bf V$ whose matrix in the canonical basis is $x$,
which is an irreducible rational representation over $\CC$.  To
simplify the notation, we denote $\rho(g)v=gv$ for every $g \in {\bf
  G}$ and $v\in{\bf V}$.  Let $\bf P$ be the stabiliser in $\bf G$ of
the line generated by $e_1$, which is a maximal (proper) parabolic
subgroup of $\bf G$ defined over $\QQ$. With $I_{k}$ the identity
$k\times k$ matrix and $s\in\RR$, let
${\bf U}=\Big\{\Big(\!\!\begin{array}{cc} 1 & u\\
  0 & I_{n-1}\end{array}\!\!\Big)\;:u\in\M_{1,n-1}(\CC)\Big\}$,
$a_s=\Big(\!\!\begin{array}{cc} e^{\frac{s}{n}}& 0\\
  0 & e^{-\frac{s}{n(n-1)}}I_{n-1}\end{array}\!\!\Big)$, and ${\bf M}=
\Big\{\Big(\!\!\begin{array}{cc} 1 & 0\\
  0 & m\end{array}\!\!\Big)\;:\;m\in\SL_{n-1}(\CC)\Big\}$. With ${\bf
  A}$ the centraliser of $\bf M$ in $\bf G$, we have that ${\bf P} =
{\bf AMU}$ is a relative Langlands decomposition of $\bf P$ over
$\QQ$, and the identity component of ${\bf A} (\RR)$ is the
one-parameter subgroup $(a_s)_{s\in\RR}$. With $\UUU$ the Lie algebra
of ${\bf U}(\RR)$, an immediate computation gives
\begin{equation}\label{eq:computlambdaappliun}
\lambda=\log \det(\Ad a_1)_{\mid \UUU}= 1>0\;.
\end{equation}

Since $F$ is homogeneous, as $F$ takes a positive value (and
equivalently), there exists $v_0\in\ZZ^n$ such that $F(v_0)>0$. We may
assume that $v_0$ is primitive up to rescaling it, and after an
integral linear change of variable (which does not change the set of
integral representations of a real number by $F$), we may assume that
$v_0=e_1$. Note that the stabiliser of $v_0$ in $\bf G$ is then
precisely $\bf MU$.

We denote by $\bf L$ the stabiliser of $F$ in $\bf G$ and by $\pi:{\bf
  L}\ra \GL({\bf V})$ the restriction of $\rho$ to $\bf L$. By the
linear independence over $\CC$ assumption, $\bf L$ is a maximal
algebraic torus defined over $\QQ$ in ${\bf G}$ (hence $\bf L$ is
reductive, but not semisimple).  For every $z\in\CC-\{0\}$, the group
$\bf L$ acts simply transitively on the affine hypersurface
$F^{-1}(z)$. Hence, with $v_s=a_sv_0=e^{\frac{s}{n}}v_0$, the orbit
\begin{equation}\label{eq:relatXsFmk}
{\bf X}_s={\bf L}v_s=F^{-1}(F(v_s))=F^{-1}(e^sF(v_0))
\end{equation}
(since $F$ is homogeneous of degree $n$) is Zariski-closed in $\bf V$.

\medskip Let us now check in two lemmas that the hypotheses of Theorem
\ref{theo:intromaincourtsc} are satisfied by these objects.

\blemm\label{lem:anisotropenormform}
The algebraic torus $\bf L$ is anisotropic over $\QQ$.
\elemm

\dem By \cite[Theo.~2.3.3, page 38]{Koch00}, there exist $a\in
\QQ-\{0\}$ and linearly independant linear forms $\ell_1,\dots,\ell_n$
on $\CC^n$ with coefficients in $\overline{\QQ}$ such that
$F=a\prod_{i=1}^n\ell_i$ and the absolute Galois group $\GalQ$ acts
transitively on the set $\{\ell_1,\dots,\ell_n\}$.  Let $\B$ be the
basis of $\CC^n$ whose dual basis is $(\ell_1,\dots,\ell_n)$. The
algebraic torus $\bf L$ is the subgroup of the elements of $\bf G$
whose matrix in the basis $\B$ is diagonal. For $1\leq i\leq n$, let
$\chi_i$ be the character (defined over $\overline{\QQ}$) of $\bf L$
which associates to an element of ${\bf L}$ the $i$-th diagonal
element of its matrix in $\B$. Note that $\GalQ$ acts transitively on
the set $\{\chi_1,\dots,\chi_n\}$. Any character of $\bf L$ may be
uniquely written $\prod_{i=1}^n\chi_i^{k_i}$ with $k_1,\dots,
k_n\in\ZZ$. Any $\QQ$-character $\prod_{i=1}^n\chi_i^{k_i}$ of $\bf
L$, being invariant under $\GalQ$, should have $k_1=\dots=k_n$ by
transitivity, hence is trivial. The result follows, since an algebraic
torus defined over $\QQ$ without nontrivial $\QQ$-characters is
anisotropic over $\QQ$, that is, it contains no nontrivial $\QQ$-split
torus (see for instance \cite[page 121]{Borel91}, though this
reference uses a different meaning of anisotropic).  \cqfd

\blemm\label{lem:Zariskiopennormform} 
The intersection ${\bf L}\cap {\bf P}$ is finite and $\bf LP$ is
Zariski-open in $\bf G$.  
\elemm

\dem Let us prove that the algebraic group ${\bf L}\cap {\bf P}$ is
finite.  Since an algebraic group has only finitely many components,
we only have to prove that its identity component ${\bf S}=({\bf L}
\cap {\bf P})_0$ is trivial.  The algebraic torus ${\bf S}$ is defined
over $\QQ$, hence is contained in a maximal torus of $\bf P$ defined
over $\QQ$.  By \cite[Theo.~19.2]{Borel91}, two maximal tori of $\bf
P$ defined over $\QQ$ are conjugated over $\QQ$.  Since $\bf G$ splits
over $\QQ$, this implies that ${\bf L}\cap {\bf P}$ splits over
$\QQ$. Since $\bf L$ is anisotropic over $\QQ$ by Lemma
\ref{lem:anisotropenormform}, this implies that ${\bf S}$ is trivial,
and proves the first claim.

Now, the homogeneous space ${\bf G}/{\bf P}$ is identified with the
complex projective space $\PP(\CC^n)$ by the map $g\mapsto \CC ge_1$.
We write $e_1=\sum_{i=1}^n c_iw_i$ where $\B=(w_i)_{1\leq i \leq n}$
is a diagonalisation basis of $\bf V$ for the action of the algebraic
torus $\bf L$, as in the proof of Lemma \ref{lem:anisotropenormform}.
Since the Galois group $\GalQ$ acts transitively on $\{w_1,\dots
w_n\}$ and fixes $e_1$, it follows that the coefficients $c_i$ are all
different from $0$. Hence ${\bf L} (\CC e_1)=\{\CC\sum_{i=1}^n b_iw_i:
b_i\ne 0\}$, which implies the second claim.  \cqfd

\medskip To conclude the proof of Corollary \ref{theo:intronormform}, we
relate the two counting functions in the statements of Corollary
\ref{theo:intronormform} and Theorem \ref{theo:intromaincourtsc}.

For every $s>0$ and $p\in\NN-\{0\}$, let $A_s^{(p)}$ be the set of
integral points of ${\bf X}_s$ whose coefficients have their greatest
common divisor equal to $p$. Note that $A_s^{(1)}={\bf X}_s\cap
\Lambda^{\rm prim}$ is the set of primitive integral points of ${\bf
  X}_s$.  With $N_s^{(p)}=\card({\bf L}(\ZZ)\bs A_s^{(p)})$, we have
$\card({\bf L}(\ZZ)\bs{\bf X}_s(\ZZ))=\sum_{p=1}^{+\infty} N_s^{(p)}$,
and $N_s^{(p)}=N_{s-\ln (p^n)}^{(1)}$, since ${\bf X}_{s-\log(p^n)}=
\frac{1}{p}{\bf X}_s$ by the homogeneity of $F$ and Equation
\eqref{eq:relatXsFmk}.

Since $\bf L$ acts simply transitively on each ${\bf X}_s$, the
stabiliser ${\bf L}_x$ of every $x\in {\bf X}_s$ is trivial, hence the
Siegel weight ${w}_{{\bf L},\pi}(x)$ is constant, equal to $\frac{1}{
  \vol( {\bf L}(\ZZ)\bs {\bf L}(\RR)}$.  By Theorem
\ref{theo:intromaincourtsc} and Equation
\eqref{eq:computlambdaappliun}, there exist $\delta>0$, that we may
assume to be in $]0,1- \frac{1}{n}[$, and $c>0$ such that, as $t\geq
0$ and $t\ra +\infty$,
$$ 
\sum_{s\in[0,t]}N_{s}^{(1)}=c\;e^{t}+\operatorname{O}(e^{t(1-\delta)})\;.
$$

For every $r\geq F(v_0)+1$, by setting $t=\log\frac{r}{F(v_0)}\geq 0$
and by using the change of variables $k=e^{s}F(v_0)$ (see Equation
\eqref{eq:relatXsFmk}), we have, with $\Sigma_k=F^{-1}(k)\cap\ZZ^n$
and $\zeta$ Riemann's zeta function,
\begin{align*}
\sum_{k\in[F(v_0),r]}\card({\bf L}(\ZZ)\bs \Sigma_k)&= 
\sum_{s\in[0,t]}\card({\bf
    L}(\ZZ)\bs{\bf X}_s(\ZZ)) =  \sum_{s\in[0,t]}\;\sum_{p=1}^{+\infty}
  N_s^{(p)}\\ &= \sum_{p=1}^{+\infty} \sum_{s\in[0,t]} N_{s-\ln
    (p^n)}^{(1)}=\sum_{p=1}^{+\infty}
  c\;p^{-n}\;e^{t}+\operatorname{O}(p^{n(\delta-1)}e^{t(1-\delta)})\\ &=
c\;\zeta(n)\;e^{t}+\operatorname{O}\big(e^{t(1-\delta)}\big)=
\frac{c\;\zeta(n)}{F(v_0)}\; r+\operatorname{O}\big(r^{1-\delta}\big)\;.
\end{align*}
Note that  $\sum_{k\in[\min\{1,F(v_0)\},\max\{1,F(v_0)\}]}\card({\bf
  L}(\ZZ)\bs \Sigma_k)$ is finite. The result follows.  
\cqfd

\subsection{\!Counting inequivalent integral points on hyperplane
sections of affine quadratic surfaces}
\label{sect:quadratic}

Let $n\geq 3$, let $q:\CC^n\ra\CC$ with $q(x)=\sum_{i=1}^n
q_{ij}x_ix_j$ for every $x=(x_1,\dots, x_n)$ be a nondegenerate
quadratic form in $n$ variables with coefficients $q_{ij}$ in $\QQ$,
and let $\ell:\CC^n\ra\CC$ with $\ell(x)=\sum_{i=1}^n
\ell_{i}x_i$ for every $x=(x_1,\dots, x_n)$ be a nonzero
linear form in $n$ variables with coefficients $\ell_{i}$ in $\QQ$.

The aim of this section is to count the number of orbits of integral
points on the sections, by the hyperplanes parallel to the kernel of
$\ell$, of the isotropic cone $q^{-1}(0)$ of $q$.

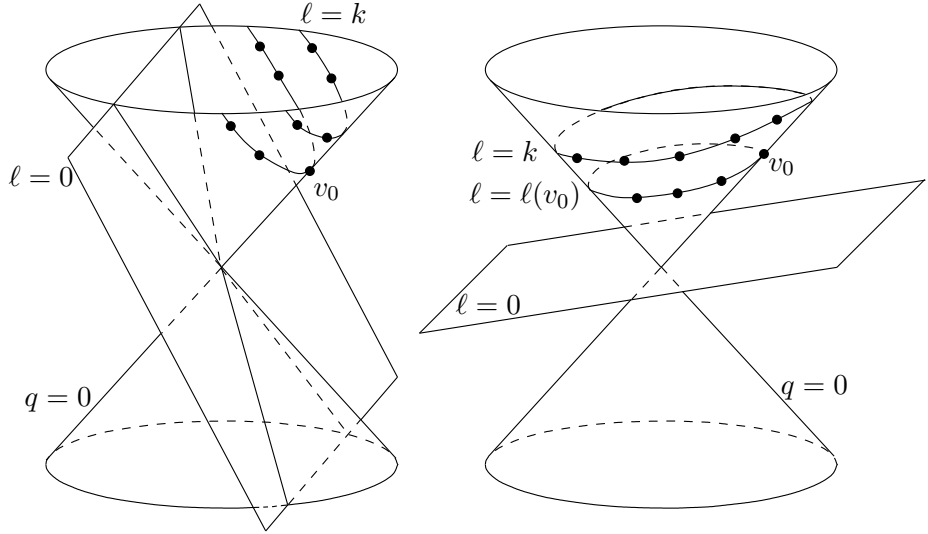
\begin{center}
\begin{picture}(0,0)%
\includegraphics{fig_quadraticsurface.pstex}%
\end{picture}%
\setlength{\unitlength}{3631sp}%
\begingroup\makeatletter\ifx\SetFigFont\undefined%
\gdef\SetFigFont#1#2#3#4#5{%
  \reset@font\fontsize{#1}{#2pt}%
  \fontfamily{#3}\fontseries{#4}\fontshape{#5}%
  \selectfont}%
\fi\endgroup%
\begin{picture}(6286,3629)(-119,-3225)
\put(3057,-683){\makebox(0,0)[lb]{\smash{{\SetFigFont{11}{13.2}{\rmdefault}{\mddefault}{\updefault}{\color[rgb]{0,0,0}$\ell=k$}%
}}}}
\put(5084,-753){\makebox(0,0)[lb]{\smash{{\SetFigFont{11}{13.2}{\rmdefault}{\mddefault}{\updefault}{\color[rgb]{0,0,0}$v_0$}%
}}}}
\put(5168,-2278){\makebox(0,0)[lb]{\smash{{\SetFigFont{11}{13.2}{\rmdefault}{\mddefault}{\updefault}{\color[rgb]{0,0,0}$q=0$}%
}}}}
\put(1976,-922){\makebox(0,0)[lb]{\smash{{\SetFigFont{11}{13.2}{\rmdefault}{\mddefault}{\updefault}{\color[rgb]{0,0,0}$v_0$}%
}}}}
\put(1887,257){\makebox(0,0)[lb]{\smash{{\SetFigFont{11}{13.2}{\rmdefault}{\mddefault}{\updefault}{\color[rgb]{0,0,0}$\ell=k$}%
}}}}
\put(3041,-972){\makebox(0,0)[lb]{\smash{{\SetFigFont{11}{13.2}{\rmdefault}{\mddefault}{\updefault}{\color[rgb]{0,0,0}$\ell=\ell(v_0)$}%
}}}}
\put(2949,-1730){\makebox(0,0)[lb]{\smash{{\SetFigFont{11}{13.2}{\rmdefault}{\mddefault}{\updefault}{\color[rgb]{0,0,0}$\ell=0$}%
}}}}
\put(-104,-852){\makebox(0,0)[lb]{\smash{{\SetFigFont{11}{13.2}{\rmdefault}{\mddefault}{\updefault}{\color[rgb]{0,0,0}$\ell=0$}%
}}}}
\put(  4,-2360){\makebox(0,0)[lb]{\smash{{\SetFigFont{11}{13.2}{\rmdefault}{\mddefault}{\updefault}{\color[rgb]{0,0,0}$q=0$}%
}}}}
\end{picture}%

\end{center}

For $\KK=\RR$ or $\QQ$, recall that $q$ is {\it isotropic (or
  indefinite when $\KK=\RR$) over $\KK$} or {\it represents $0$ over
  $\KK$} if there exists $x\in \KK^n-\{0\}$ such that $q(x)=0$, and
that $q$ is {\it anisotropic over $\KK$} otherwise. For instance,
$x^2+2y^2-7z^2$ is anisotropic over $\QQ$, but indefinite over
$\RR$. By A.~Meyer's 1884 result (see for instance \cite[page
77]{Serre70}), if $n\geq 5$, then $q$ is isotropic over $\QQ$ if and
only if $q$ is indefinite over $\RR$.

\bigskip
\noindent{\bf Proof of Corollary \ref{theo:introquadraticsurface}. } 
In order to apply Theorem \ref{theo:mainmaximal} (2), let us first
define the objects appearing in its statement.

Let ${\bf G}=O_q$ be the orthogonal group of the nondegenerate
rational quadratic form $q$, which is a connected semisimple linear
algebraic group defined over $\QQ$, hence is reductive without
nontrivial $\QQ$-characters. Let ${\bf V}=\CC^n$ and let $\rho:{\bf
  G}\ra \GL({\bf V})$ be the monomorphism mapping a matrix $x$ to the
linear automorphism of $\bf V$ whose matrix in the canonical basis is
$x$, which is an irreducible rational representation over $\CC$. Let
$\Lambda=\ZZ^n$, which is a $\ZZ$-lattice in ${\bf V}(\QQ)$ invariant
under ${\bf G}(\ZZ)$. To simplify the notation, we denote
$\rho(g)v=gv$ for every $g \in {\bf G}$ and $v\in{\bf V}$.

Since $q$ is assumed to be isotropic over $\QQ$, there exists $v_0$ in
$\Lambda-\{0\}$ such that $q(v_0)=0$ and we assume that $\ell(v_0)\geq
0$ up to replacing $v_0$ by $-v_0$.  Since the restriction of $q$ to
the kernel of $\ell$ is assumed to be anisotropic over $\QQ$, we have
$\ell(v_0)>0$. Let $\bf P$ be the stabiliser in $\bf G$ of the line
generated by $v_0$, which is a maximal (proper) parabolic subgroup of
$\bf G$ defined over $\QQ$ since this line is isotropic. Let
$\B=(e_1,\dots, e_n)$ be a basis of ${\bf V}$ over $\QQ$ such that
$e_1=v_0$, $(e_1,e_2)$ is a standard basis of a hyperbolic plane over
$\QQ$ for $q$, which is orthogonal for $q$ to the vector subspace
${\bf V}'$ generated by $\B'=(e_3,\dots, e_n)$. In particular, the
matrix of $q$ in the basis $\B$ is $Q= \left(\!\!\begin{array}{ccc}
    0&1& 0\\ 1 & 0 & 0 \\ 0 & 0 & Q'
  \end{array}\!\!\right)$ with $Q'$ the (rational symmetric) matrix in
the basis $\B'$ of the restriction $q'$ of $q$ to ${\bf V}'$. Denoting
in the same way a vector $v$ (resp.~$u$) of ${\bf V}$ (resp.~${\bf
  V}'$) and the column vector of its coordinates in $\B$
(resp.~$\B'$), we have $q(v)={\,}^tvQv$ (resp.~$q'(u)={\,}^tuQ'u$).
With $I_{k}$ the identity $k\times k$ matrix and $s\in\RR$, define
$$
a_s=\left(\!\!\begin{array}{ccc} 
e^s&0& 0\\ 0& e^{-s}  & 0 \\ 0 & 0 & I_{n-2}
\end{array}\!\!\right)\;,\;\;\;{\bf A}=
\left\{\left(\!\!\begin{array}{ccc} 
a&0& 0\\ 0 &a^{-1}  & 0 \\ 0 & 0 & I_{n-2}
\end{array}\!\!\right)\;:\;a\in\CC^*\right\}\;,$$ 
$${\bf M}= \left\{\left(\!\!\begin{array}{ccc} 
1&0& 0\\ 0 &1  & 0 \\ 0 & 0 & m
\end{array}\!\!\right)\;:\;m\in O_{q'}\right\}\;\;\;
 {\rm and}\;\;\;
{\bf U}=\left\{\left(\!\!\begin{array}{ccc} 
1&-q'(u)/2& -{\,}^tuQ'\\ 0 &1  & 0 \\ 0 & u & I_{n-2}
\end{array}\!\!\right)\;:\;u\in{\bf V}'\right\}\;.
$$
It is easy to check that ${\bf P} ={\bf AMU}$ is a relative Langlands
decomposition of $\bf P$, that the identity component of ${\bf A}
(\RR)$ is the one-parameter subgroup $(a_s)_{s\in\RR}$, and that the
stabiliser of $v_0=e_1$ in $\bf G$ is exactly ${\bf MU}$. With $\UUU$
the Lie algebra of ${\bf U}(\RR)$, an immediate computation gives
(since $n\geq 3$)
\begin{equation}\label{eq:computlambdaapplideu}
\lambda=\log \det(\Ad a_1)_{\mid \UUU}= n-2>0\;.
\end{equation}

We denote by ${\bf L}=\{g\in{\bf G}\;:\;\ell\circ g = \ell\}$ the
stabiliser of $\ell$ in $\bf G$, which is a linear algebraic group
defined over $\QQ$. Let ${\bf W}$ be the kernel of $\ell$ and ${\bf
  W}^\perp$ be its orthogonal for $q$. Since $q_{\mid {\bf W}}$ is
assumed to be nondegenerate, ${\bf W}^\perp$ is a line, ${\bf V} =
{\bf W}^\perp\oplus{\bf W}$, and the bloc matrix of $q$ in this
decomposition is diagonal.

\medskip Let us now check in the next lemma that the hypotheses of
Theorem \ref{theo:mainmaximal} are satisfied by these objects.

\blemm\label{lem:anisotropequadsurf} (1) The linear algebraic group $\bf
L$ is reductive and anisotropic over $\QQ$.  

(2) For every $s\in\RR$, if $k=e^s\ell(v_0)$ and ${\bf X}_s={\bf L}
a_sv_0$, then ${\bf X}_s=\{v\in{\bf V}\;:\;q(v)=0,\ell(v)=k\}$. In
particular, ${\bf X}_s$ is Zariski-closed in ${\bf V}$. 

(3) The subset $\bf LP$ is Zariski-open in $\bf G$.  \elemm

\dem (1) For every $g\in\GL({\bf V})$, if $\ell\circ g=\ell$, then $g$
preserves ${\bf W}$. If furthermore $g\in{\bf G}=O_q$, then $g$
preserves ${\bf W}^\perp$. Since ${\bf W}^\perp$ is a line, there
exists $\lambda\in\CC$ such that $g$ acts by $x\mapsto \lambda x$ on
${\bf W}^\perp$. As $\ell_{\mid {\bf W}^\perp}$ is nonzero and $g$
preserves $\ell$, we have $\lambda=1$. Hence the elements of $\bf L$
are exactly the elements of $\GL({\bf V})$ whose bloc matrix in the
decomposition ${\bf V} ={\bf W}^\perp\oplus{\bf W}$ has the form
$\Big(\!\!\begin{array}{cc} 1 & 0\\ 0 & g'\end{array} \!\!\Big)$ 
with $g'\in O_{q_{\mid {\bf W}}}$.  In particular, the linear
algebraic group $\bf L$, isomorphic over $\QQ$ to the orthogonal group
of the nondegenerate rational quadratic form $q_{\mid {\bf W}}$, is
semisimple hence reductive.

It is well-known (see for instance \cite{Borel62}\cite[page
270]{BorJi06}) that the $\QQ$-rank of the orthogonal group $O_{q''}$
of a nondegenerate rational quadratic form $q''$ is zero (or
equivalently that $O_{q''}$ is anisotropic over $\QQ$) if and only if
$q''$ does not represents $0$ over $\QQ$. For instance, this follows
from the fact that the spherical Tits building over $\QQ$ of $O_{q''}$
is the building of isotropic flags over $\QQ$. Hence by assumption,
$\bf L$ is anisotropic over $\QQ$.

\medskip (2) Note that by the definition of $a_s$, we have
$a_sv_0=e^sv_0$, hence by the linearity of $\ell$, we may assume that
$s=0$. Recall that $q(v_0)=0$ and $\ell(v_0)> 0$.  By the definition
of $\bf L$, the orbit ${\bf X}_0={\bf L}v_0$ is contained in
$\{v\in{\bf V}\;:\; q(v)=0,\; \ell(v)=\ell(v_0)\}$. To prove the
opposite inclusion, write $v=v'+v''$ the decomposition of any
$v\in{\bf V}$ in the direct sum ${\bf V} ={\bf W}^\perp\oplus{\bf
  W}$. If $\ell(v)=\ell(v_0)$ and $q(v)=0$, then $v'=v'_0$ and
$q(v'')=-q(v')=-q(v'_0)$, and in particular $q(v'')=q(v_0'')$. By
Witt's theorem, there exists $g'\in O_{ q_{\mid {\bf W}}}$ such that
$v''=g'v_0''$. Hence the linear transformation of $\bf V$ which is
the identity on ${\bf W}^\perp$ and is equal to $g'$ on ${\bf W}$, is
an element of $\bf L$ sending $v=v'+v''$ to $v_0=v'_0+v_0''$. The
second assertion follows.

\medskip (3) The algebraic group ${\bf G}=O_q$ acts transitively on
the projective variety of isotropic lines in $\bf V$, the stabiliser
of the line generated by $v_0$ being $\bf P$ by definition.  As we
have seen in (2), the orbit under $\bf L$ of the line
generated by $v_0$ is hence the Zariski-open subset of ${\bf G}/{\bf
  P}$ consisting of the isotropic lines not contained in ${\bf
  W}$. The last claim of Lemma \ref{lem:anisotropequadsurf} follows.
\cqfd

\medskip To conclude the proof of Corollary
\ref{theo:introquadraticsurface}, we relate the two counting functions
in the statements of Corollary \ref{theo:introquadraticsurface} and
Theorem \ref{theo:mainmaximal} (2). Let $L={\bf L}(\RR)_0$ and
$\Ga={\bf G}(\RR)_0\cap{\bf G}(\ZZ)$.

We have $\ell(v_0)> 0$ by the definition of $v_0$.  For every $r\geq
\ell(v_0)+1$, let $t=\ln \frac{r}{\ell(v_0)} >0$. With $\Sigma_k$ as
in the statement of Corollary \ref{theo:introquadraticsurface}, using
the change of variables $k= e^s\ell(v_0)$ and Lemma
\ref{lem:anisotropequadsurf} (2), by the definition of the modified
Siegel weights in Equation \eqref{eq:defimodifSiegelweight}, we have
$$
\sum_{k\in[\ell(v_0),r]}\;\sum_{[u]\in ({\bf L}(\ZZ)\cap L)\bs\Sigma_k}
\vol\big( ({\bf L}_u(\ZZ)\cap L)\bs({\bf L}_u\cap L)\big) =
$$
\begin{equation}\label{eq:relatcoubntingquadsurf}
\vol\big( (L\cap\Ga)\bs L\big)\sum_{s\in[0,t]}
\;\sum_{[u]\in ({\bf L}(\ZZ)\cap L)\bs{\bf X_s}\cap\Lambda^{\rm prim}}
{w'}_{{\bf L},\rho_{\mid {\bf L}}}(u)\;.
\end{equation}
By Theorem \ref{theo:mainmaximal} (2) and Equation
\eqref{eq:computlambdaapplideu}, there exist $c,\delta>0$ such that as
$t\ra+\infty$, the quantity \eqref{eq:relatcoubntingquadsurf} is equal
to
$$
c\;e^{(n-2)t}+\operatorname{O}\big(e^{(n-2-\delta)t}\big)=
\frac{c}{\ell(v_0)^{n-2}}\;r^{n-2}+
\operatorname{O}\big(r^{n-2-\delta}\big)\;.
$$
Note that $\sum_{k\in[\min\{1,\ell(v_0)\},\max\{1,\ell(v_0)\}]}\;
\sum_{[u]\in ({\bf L}(\ZZ)\cap L) \bs\Sigma_k} \vol\big( ({\bf L}_u
(\ZZ) \cap L)\bs({\bf L}_u\cap L)\big)$ is finite. This concludes the
proof of Corollary \ref{theo:introquadraticsurface}. \cqfd

\medskip \noindent {\bf Remarks } (1) If $n\geq 6$, since $q$ is
isotropic over $\QQ$ and the restriction of $q$ to the kernel of
$\ell$ is anisotropic over $\QQ$, then the signature of $q$ over $\RR$
is $(1,n-1)$ or $(n-1,1)$, and ${\bf L} (\RR)$ is compact (see the
above picture on the right); hence ${\bf L} (\ZZ)$ is finite, and our
result allows to count integral points on the quadratic hypersurface
$q^{-1}(0)$ (see the references given in the introduction for related
works).

\medskip (2) If $n\geq 4$, then we have a result similar to
Corollary \ref{theo:introquadraticsurface} where we consider all the
integral points and not only the primitive ones: under the other
assumptions of Corollary \ref{theo:introquadraticsurface} and with $c$
as above, we have, for every $r\geq 1$
with $r\ra+\infty$,
\begin{align*}
\sum_{k\in[1,r]} \;\;\;\sum_{[u]\in ({\bf L}(\ZZ)\cap L)\bs
  (q^{-1}(0)\cap \ell^{-1}(k)\cap\ZZ^n)}\; &\vol\big( ({\bf
  L}_u(\ZZ)\cap L)\bs({\bf L}_u\cap L)\big) \\ & 
=\;\; \frac{c\;\zeta(n-2)}{\ell(v_0)^{n-2}} \,r^{n-2}+
\operatorname{O}\big(r^{n-2-\delta}\big)\;.
\end{align*}
The proof is similar to the one at the end of Section
\ref{sect:normform}.  For every $s\in\RR$ and $p\in\NN-\{0\}$, we
denote by $A_s^{(p)}$ the set of integral points of ${\bf X}_s$ whose
greatest common divisor of their coefficients is $p$. We note that by
Lemma \ref{lem:anisotropequadsurf} (2), the map from $A_s^{(p)}$ to
$A_{s-\ln p}^{(1)}$ defined by $x\mapsto \frac{x}{p}$ is a bijection
such that ${\bf L}_{\frac{x}{p}}={\bf L}_x$ for every $x\in
A_s^{(p)}$. Hence with
$$
N_s^{(p)}=\sum_{[u]\in ({\bf L}(\ZZ)\cap L)\bs A_s^{(p)}}\; 
\vol\big( ({\bf L}_u(\ZZ)\cap L) \bs ({\bf L}_u (\RR)\cap L)\big)\;,
$$ 
we have $N_s^{(p)}= N_{s-\ln p}^{(1)}$ and
$$
\sum_{[u]\in ({\bf L}(\ZZ)\cap L)\bs (q^{-1}(0)\cap
  \ell^{-1}(k)\cap\ZZ^n)}\vol\big( ({\bf L}_u(\ZZ)\cap L) \bs ({\bf L}_u
\cap L)\big)=\sum_{p=1}^{\infty}N_s^{(p)}\;,
$$
and one concludes as in the end of Section \ref{sect:normform}.

When $n=3$, the same argument gives
$$
\sum_{k\in[1,r]} \;\;\;\sum_{[u]\in ({\bf L}(\ZZ)\cap L)\bs
  (q^{-1}(0)\cap \ell^{-1}(k)\cap\ZZ^n)}\; \vol\big( ({\bf
  L}_u(\ZZ)\cap L)\bs({\bf L}_u\cap L)\big) = \frac{c}{\ell(v_0)}\;
\,r\log r+ \operatorname{O}\big(r\big)\;.
$$

\subsection{Counting inequivalent integral points of 
given norm in central division algebras}
\label{sect:unitdivisionalgebra}

Let $n\geq 2$, let $D$ be a central simple algebra over $\QQ$ of
dimension $n^2$, let $N:D\ra \QQ$ be its reduced norm, and let $\OOO$
be an {\it order} in $D$ (that is, a finitely generated
$\ZZ$-submodule of $D$, generating $D$ as a $\QQ$-vector space, which
is a unitary subring). We refer for instance to \cite{Reiner75} and
\cite[Chap.~I, \S 1.4]{PlaRap94}) for generalities.  The aim of this
section is to use our main result to deduce asymptotic counting
results of elements of $\OOO$ (modulo units) of given norm.  We note
that this result can be also deduced by elementary methods (as for the
first example) of fundamental domain for the action of unit groups on
a level set of the norm, as well as from the analytic continuation
(established in \cite{SatShi74}) of the zeta function associated to
the corresponding prehomogeneous vector space.

\bcoro\label{theo:applitroi} 
If $D$ is a division algebra over $\QQ$, then there exist
$c=c(D,\OOO)>0$ and $\delta=\delta(D)>0$ such
that, for every $r\ge 1$ with $r\ra+\infty$,
$$
\card\;_{\OOO^\times\bs}\{x\in\OOO\;:\;1\leq  |N(x)|\leq r\}= 
c\; r^{n}\big(1+
\operatorname{O}\big(r^{-\delta}\big)\big)\;.
$$
\ecoro

\dem In order to apply Theorem \ref{theo:intromaincourtsc}, let us
first define the objects appearing in its statement.

Let $\bf V$ be the vector space over $\QQ$ such that ${\bf V}(K) =
D\otimes_\QQ K$ for every characteristic zero field, with the integral
structure such that $\Lambda={\bf V}(\ZZ)=\OOO$, which is (for the
extended multiplication) a central simple algebra over $\CC$. Let ${\bf
  D}^1$ be the group of elements of (reduced) norm $\pm 1$ in $\bf V$.

We take ${\bf G}=\SL({\bf V})$ (which is connected, simply connected,
semisimple, defined over $\QQ$, hence reductive without nontrivial
$\QQ$-characters) and $\rho$ the inclusion of ${\bf G}$ in
$\operatorname{GL}({\bf V})$ (which is an irreducible rational
representation).  To simplify the notation, we denote $\rho(g)v=gv$
for every $g \in {\bf G}$ and $v\in{\bf V}$.

Let ${\bf L}$ be the algebraic subgroup of ${\bf G}$ which is the
image of ${\bf D}^1$ into ${\bf G}$ by the (left) regular
representation $d\mapsto \{v\mapsto dv\}$. Note that the linear
algebraic groups ${\bf L}$ and ${\bf D}^1$ are defined over $\QQ$ and
are isomorphic by this representation. We have
\begin{equation}\label{eq:computinversOOO}
{\bf L}(\ZZ)={\bf D}^1\cap \OOO=\OOO^\times\;.
\end{equation}

We take $v_0\in {\bf V}$ to be the identity element in $D$. The
stabiliser of the line $\CC v_0$ in ${\bf G}$ is a (maximal) parabolic
subgroup ${\bf P}$ of ${\bf G}$ defined over $\QQ$.  We note that
$\dim({\bf P})=\dim(D)^2-\dim(D)-1$ and $\dim({\bf L})=\dim(D)-1$. We
have a relative Langlands decomposition ${\bf P}={\bf AMU}$ with ${\bf
  MU}$ the stabiliser of $v_0$ in $\bf G$, and we may write ${\bf A}
(\RR)_0=(a_s)_{s\in\RR}$ such that $a_sv_0=e^{\frac{s}{n}}v_0$. An easy 
computation gives
\begin{equation}\label{eq:computlambdaapplitroi}
\lambda=\log \det(\Ad a_1)_{\mid \UUU}= n>0\;.
\end{equation}

\medskip Let us now check that the hypotheses of Theorem
\ref{theo:intromaincourtsc} are satisfied by these objects.

We claim that the group ${\bf L}\cap {\bf P}$ is finite.  The action
of this group on $v_0$ defines a $\mathbb{Q}$-character of ${\bf
  L}\cap {\bf P}$. Since ${\bf L}\simeq {\bf D}^1$ is anisotropic over
$\QQ$ (see for instance \cite[Chap.~II, \S 2.3]{PlaRap94}), this
character must be trivial on $({\bf L}\cap {\bf P})_0$, and $({\bf
  L}\cap {\bf P})_0v_0=v_0$. Since $\hbox{Stab}_{\bf L}(v_0)=\{e\}$,
it follows that $({\bf L}\cap {\bf P})_0=\{e\}$, which proves the
claim. Comparing dimensions, we deduce that ${\bf L} {\bf P}$ is
Zariski-open in ${\bf G}$.

For every $s\in\RR$, we have
\begin{equation}\label{eq:coputXsapplitroi}
{\bf X}_s={\bf L}a_s v_0=e^{\frac{s}{n}}{\bf L}v_0=
e^{\frac{s}{n}}{\bf D}^1\;.
\end{equation}
Hence ${\bf X}_s$ is Zariski-closed in $\bf V$. 

\medskip
To conclude the proof of Corollary \ref{theo:applitroi}, we
relate the two counting functions in the statements of Corollary
\ref{theo:applitroi} and Theorem \ref{theo:intromaincourtsc}.

Since $\bf L$ acts simply transitively on the orbit of $v_0$, the
Siegel weights are constant, equal to $\frac{1}{\vol({\bf L}(\ZZ)\bs
  {\bf L}(\RR))}$. For every $k\in\NN-\{0\}$, denote by $\OOO^{(k)}$
the subset of nonzero elements of $\OOO$ whose greatest common divisor
of their coefficients in a $\ZZ$-basis of $\OOO$ is $k$. In
particular, since the norm is a homogeneous polynomial of degree $n$
and by Equation \eqref{eq:coputXsapplitroi}, we have
$$
{\bf X}_s\cap \Lambda^{\rm prim}=\{x\in \OOO^{(1)}\;:\; N(x)=e^s\}\;.
$$
Note that the map $x\mapsto \frac{x}{k}$ is a bijection from
$\OOO^{(k)}$ to $\OOO^{(1)}$. Hence, using Equation
\eqref{eq:computinversOOO} and Theorem \ref{theo:intromaincourtsc},
there exist $\delta>0$, that we may assume to be in $]0,1[\,$, and
$c>0$ such that, as $r\geq 1$ and $r\ra +\infty$,
\begin{align*}
  \card\;_{\OOO^\times\bs}\{x\in\OOO\;:\;1\leq |N(x)|\leq r\} & =
  \sum_{k=1}^{+\infty}\card\;_{\OOO^\times\bs}\{x\in\OOO^{(k)}\;:\;1\leq
  |N(x)|\leq r\} \\ & =
  \sum_{k=1}^{+\infty}\card\;_{\OOO^\times\bs}\{x\in\OOO^{(1)}\;:\;1\leq
  |N(x)|\leq \frac{r}{k^n}\}\\ & =
  \sum_{k=1}^{+\infty}\;\;\sum_{0\leq s\leq \log\frac{r}{k^n}}\;
\card\;\big({\bf L}(\ZZ)\bs\big({\bf X}_s\cap\Lambda^{\rm
  prim}\big)\big)
\\ & =
  \sum_{k=1}^{+\infty}\;c\;\;\big(\frac{r}{k^n}\big)^{n}
\Big(1+\operatorname{O}\big(\big(\frac{r}{k^n}\big)^{-\delta}\big)\Big)
\\ & =
 c\;\zeta(n^2)\;r^{n}
\big(1+\operatorname{O}(r^{-\delta})\big)\;.
\end{align*}
This ends the proof of Corollary \ref{theo:applitroi}.
\cqfd

{\small \bibliography{../biblio} }

\begin{thebibliography}{BdlHV}

\bibitem[Bab]{Babillot02a}
M.~Babillot.
\newblock {\it Points entiers et groupes discrets~: de l'analyse aux syst\`emes
  dynamiques}.
\newblock {in "Rigidit\'e, groupe fondamental et dynamique", Panor. Synth\`eses
  {\bf 13}, 1--119, Soc. Math. France, 2002}.

\bibitem[BdlHV]{BekHarVal08}
B.~Bekka, P.~de~la Harpe, and A.~Valette.
\newblock {\it Kazhdan's property T}.
\newblock {New Math. Mono. {\bf 11}, Cambridge Univ. Press, 2008}.

\bibitem[BO]{BenOh12}
Y.~Benoist and H.~Oh.
\newblock {\it Effective equidistribution of $S$-integral points on symmetric
  varieties}.
\newblock {To appear in Annales de L'Institut Fourier}.

\bibitem[Bor1]{Borel62}
A.~Borel.
\newblock {\it Ensembles fundamentaux pour les groupes arithmétiques}.
\newblock {in ``Colloque sur la Théorie des Groupes Algébriques'', CBRM,
  Bruxelles, 1962, pp. 23--40}.

\bibitem[Bor2]{Borel69}
A.~Borel.
\newblock {\it Introduction aux groupes arithmétiques}.
\newblock {Hermann, 1969}.

\bibitem[Bor3]{Borel91}
A.~Borel.
\newblock {\it Linear algebraic groups}.
\newblock {2nd Enlarged Ed., Grad. Texts Math. {\bf 126}, Springer Verlag,
  1991}.

\bibitem[Bor4]{Borel66}
A.~Borel.
\newblock {\it Reduction theory for arithmetic groups}.
\newblock {in ``Algebraic Groups and Discontinuous Subgroups'', A.~Borel and
  G.~D.~Mostow eds, Proc. Sympos. Pure Math. (Boulder, 1965), pp. 20--25, Amer.
  Math. Soc. 1966}.

\bibitem[BHC]{BorelHarishChandra62}
A.~Borel and Harish-Chandra.
\newblock {\it Arithmetic subgroups of algebraic groups}.
\newblock {Ann. of Math. {\bf 75} (1962) 485--535}.

\bibitem[BJ]{BorJi06}
A.~Borel and L.~Ji.
\newblock {\it Compactifications of symmetric and locally symmetric spaces}.
\newblock {Birkh\"auser, 2006}.

\bibitem[BR]{BorRud95}
M.~Borovoi and Z.~Rudnick.
\newblock {\it Hardy-Littlewood varieties and semisimple groups}.
\newblock {Invent. Math. {\bf 119} (1995) 37--66}.

\bibitem[Clo]{Clozel03}
L.~Clozel.
\newblock {\it D\'emonstration de la conjecture $\tau$}.
\newblock {Invent. Math. {\bf 151} (2003) 297--328}.

\bibitem[CTX]{ColXu09}
J.-L. Colliot-Thélène and F.~Xu.
\newblock {\it Brauer-Manin obstruction for integral points of homogeneous
  spaces and representation by integral quadratic forms}.
\newblock {Compositio Math. {\bf 145} (2009) 309--363}.

\bibitem[Cow]{Cowling79}
M.~Cowling.
\newblock {\it Sur les coefficients des représentations unitaires des groupes
  de Lie simples}.
\newblock {in "Analyse harmonique sur les groupes de Lie" (Sém.
  Nancy-Strasbourg 1976–1978), II, pp.~132–178, Lect.~Notes Math. {\bf
  739}, Springer Verlag, 1979}.

\bibitem[DRS]{DukRudSar93}
W.~Duke, Z.~Rudnick, and P.~Sarnak.
\newblock {\it Density of integer points on affine homogeneous varieties}.
\newblock {Duke Math. J. {\bf 71} (1993) 143--179}.

\bibitem[EM]{EskMcMul93}
A.~Eskin and C.~McMullen.
\newblock {\it Mixing, counting, and equidistribution in Lie groups}.
\newblock {Duke Math. J. {\bf 71} (1993) 181--209}.

\bibitem[EMS]{EskMozSha96}
A.~Eskin, S.~Mozes, and N.~Shah.
\newblock {\it Unipotent flows and counting lattice points on homogeneous
  varieties}.
\newblock {Ann. of Math. {\bf 143}(1996) 253–-299}.

\bibitem[EO]{EskOh06}
A.~Eskin and H.~Oh.
\newblock {\it Representations of integers by an invariant polynomial and
  unipotent flows}.
\newblock {Duke Math. J. {\bf 135} (2006) 481--506}.

\bibitem[ERS]{EskRudSar91}
A.~Eskin, Z.~Rudnick, and P.~Sarnak.
\newblock {\it A proof of Siegel's weight formula}.
\newblock {Internat. Math. Res. Notices {\bf 5} (1991) 65--69}.

\bibitem[GO]{GanOh03}
W.~T. Gan and H.~Oh.
\newblock {\it Equidistribution of integer points on a family of homogeneous
  varieties: a problem of Linnik}.
\newblock {Compositio Math. {\bf 136} (2003) 323--352}.

\bibitem[Hir]{Hirsch76}
M.~Hirsch.
\newblock {\it Differential topology}.
\newblock {Grad.~Texts Math. {\bf 33}, Springer Verlag, 1976}.

\bibitem[KS]{KelSar09}
D.~Kelmer and P.~Sarnak.
\newblock {\it Strong spectral gaps for compact quotients of products of ${\rm
  PSL}(2,\RR)$}.
\newblock {J. Euro. Math. Soc. {\bf 11} (2009) 283--313}.

\bibitem[Kim]{Kimura03}
T.~Kimura.
\newblock {\it Introduction to prehomogeneous vector spaces}.
\newblock {Transl. Math. Mono. {\bf 215}, Amer. Math. Soc. 2003}.

\bibitem[KM1]{KleMar96}
D.~Kleinbock and G.~Margulis.
\newblock {\it Bounded orbits of nonquasiunipotent flows on homogeneous
  spaces}.
\newblock {Sinai's Moscow Seminar on Dynamical Systems, 141--172, Amer. Math.
  Soc. Transl. Ser. {\bf 171}, Amer. Math. Soc. 1996}.

\bibitem[KM2]{KleMar99}
D.~Kleinbock and G.~Margulis.
\newblock {\it Logarithm laws for flows on homogeneous spaces}.
\newblock {Invent. Math. {\bf 138} (1999) 451--494}.

\bibitem[Koc]{Koch00}
H.~Koch.
\newblock {\it Number theory : algebraic numbers and functions}.
\newblock {Grad. Stud. Math. {\bf 24}, Amer. Math. Soc. 2000}.

\bibitem[Lan]{Lang94}
S.~Lang.
\newblock {\it Algebraic number theory}.
\newblock {Grad Texts Math. 2nd ed., Springer Verlag, 1994}.

\bibitem[Nev]{Nevo05}
A.~Nevo.
\newblock {\it Exponential volume growth, maximal functions on symmetric
  spaces, and ergodic theorems for semi-simple Lie groups}.
\newblock {Erg.~Theo.~Dyn.~Syst. {\bf 25} (2005) 1257--1294}.

\bibitem[Oh1]{Oh04}
H.~Oh.
\newblock {\it Hardy-Littlewood system and representations of integers by an
  invariant polynomial}.
\newblock {Geom. Funct. Anal. {\bf 14} (2004) 791--809}.

\bibitem[Oh2]{Oh10}
H.~Oh.
\newblock {\it Orbital counting via mixing and unipotent flows}.
\newblock {In "Homogeneous flows, moduli spaces and arithmetic", M.~Einsiedler
  et al eds., Clay Math. Proc. {\bf 10}, Amer. Math. Soc. 2010, 339--375}.

\bibitem[PP]{ParPau12JMD}
J.~Parkkonen and F.~Paulin.
\newblock {\it \'Equidistribution, comptage et approximation par irrationnels
  quadratiques}.
\newblock {J. Mod. Dyn. {\bf 6} (2012) 1--40}.

\bibitem[PR]{PlaRap94}
V.~Platonov and A.~Rapinchuck.
\newblock {\it Algebraic groups and number theory}.
\newblock {Academic Press, 1994}.

\bibitem[Rag]{Raghunathan72}
M.~Raghunathan.
\newblock {\it Discrete subgroups of Lie groups}.
\newblock {Springer Verlag, 1972}.

\bibitem[Rei]{Reiner75}
I.~Reiner.
\newblock {\it Maximal orders}.
\newblock {Academic Press, 1972}.

\bibitem[SS]{SatShi74}
M.~Sato and T.~Shintani.
\newblock {\it On zeta functions associated with prehomogeneous vector spaces}.
\newblock {Ann. of Math. {\bf 100} (1974) 131--170}.

\bibitem[Ser]{Serre70}
J.-P. Serre.
\newblock {\it Cours d'arithmetique}.
\newblock {Press. Univ. France, Paris, 1970}.

\bibitem[Sie1]{Siegel44}
C.~L. Siegel.
\newblock {\it On the theory of indefinite quadratic forms}.
\newblock {Ann. of Math. {\bf 45} (1944) 577--622}.

\bibitem[Sie2]{Siegel44b}
C.~L. Siegel.
\newblock {\it The average measure of quadratic forms with given determinant
  and signature}.
\newblock {Ann. of Math. {\bf 45} (1944) 667--685}.

\bibitem[Spr]{Springer94}
T.~A. Springer.
\newblock {\it Linear algebraic groups}.
\newblock {In "Algebraic geometry IV", A.~Parshin, I.~Shavarevich eds., Encyc.
  math. Scien. {\bf 55}, Springer Verlag, 1994}.

\bibitem[Vos]{Voskresenskii98}
V.~E. Voskresenskii.
\newblock {\it Algebraic groups and their birational invariants}.
\newblock {Transl. Math. Mono. {\bf 179}, Amer. Math. Soc., 1998}.

\bibitem[Wei]{Weil65}
A.~Weil.
\newblock {\it L'intégration dans les groupes topologiques et ses
  applications}.
\newblock {Hermann, 1965}.

\end{thebibliography}

\bigskip
{\small\noindent \begin{tabular}{l} 
School of Mathematics \\ 
University of Bristol \\ 
Bristol BS8 1TW, United Kingdom\\
{\it e-mail: a.gorodnik@bristol.ac.uk}
\end{tabular}
\medskip

\noindent \begin{tabular}{l}
D\'epartement de math\'ematique, UMR 8628 CNRS, B\^at.~425\\
Universit\'e Paris-Sud 11,
91405 ORSAY Cedex, FRANCE\\
{\it e-mail: frederic.paulin@math.u-psud.fr}
\end{tabular}}

\end{document}